\documentclass[11pt, leqno]{amsart}
\usepackage{amsfonts,amssymb}

\usepackage{amsmath}
\usepackage{amsthm}
\usepackage{amsrefs}
\usepackage{qsymbols}
\usepackage{latexsym}
\usepackage{chngcntr}

\usepackage{paralist}

\usepackage{mathtools}

\usepackage{esint}

\usepackage[hidelinks]{hyperref}
\usepackage{comment}
\usepackage{tikz-cd}
\usepackage{import}
\usepackage{todonotes}
\usepackage{csquotes}
\usepackage[dvipsnames]{xcolor}

\usepackage{enumitem}

\newcommand*{\rom}[1]{\expandafter\@slowromancap\romannumeral #1@}

\newtheorem{theorem}{Theorem}[section]

\newtheorem{lemma}[theorem]{Lemma}
\newtheorem{proposition}[theorem]{Proposition}
\newtheorem{corollary}[theorem]{Corollary}
\theoremstyle{definition}
\newtheorem{definition}[theorem]{Definition}
\newtheorem{remark}[theorem]{Remark}

\newtheorem{assumption}[theorem]{Assumption}

\newcommand{\IR}{\mathbb{R}}
\newcommand{\IC}{\mathbb{C}}
\newcommand{\IN}{\mathbb{N}}
\newcommand{\IZ}{\mathbb{Z}}

\renewcommand{\a}{\mathbf{a}}

\newcommand{\cO}{\mathcal{O}}

\newcommand{\cF}{\mathcal{F}}

\newcommand{\cR}{\mathcal{R}}
\newcommand{\cA}{\mathcal{A}}

\newcommand{\cL}{\mathcal{L}}
\newcommand{\cD}{\mathcal{D}}

\renewcommand{\L}{\mathrm{L}}

\renewcommand{\H}{\mathrm{H}}

\renewcommand{\S}{\mathrm{S}}

\newcommand{\C}{\mathrm{C}}
\newcommand{\W}{\mathrm{W}}

\newcommand{\I}{\mathrm{I}}

\newcommand{\Lloc}{\L_{\mathrm{loc}}}

\newcommand{\ind}{{\mathbf{1}}}

\newcommand{\NL}{\mathrm{NL}}

\renewcommand{\d}{\mathrm{d}}
\newcommand{\eps}{\varepsilon}
\newcommand{\loc}{\mathrm{loc}}
\renewcommand\Re{\operatorname{Re}}

\renewcommand{\div}{\operatorname{div}}

\DeclareMathOperator{\supp}{supp}

\DeclareMathOperator{\Id}{Id}

\numberwithin{equation}{section}

\title[The Kato square root property for fractional elliptic operators]{The Kato square root property for fractional elliptic operators in divergence form}

\author{Luca Haardt}

\address{Karlsruhe Institute of Technology, Department of Mathematics, 76131 Karlsruhe, Germany}
\email{luca.haardt@kit.edu}



\date{\today}

\thanks{The author was supported by \textit{Studienstiftung des deutschen Volkes}. }

\begin{document}

\begin{abstract}
    We establish the Kato square root property for fractional elliptic operators $L = -\operatorname{div}^\alpha(\mu \nabla^\alpha \cdot)$ of order $2\alpha \in (0,2)$ on $\mathbb{R}^d$, built from the Riesz fractional gradient and divergence, with bounded, measurable and uniformly elliptic complex coefficients $\mu$. More precisely, we show that the domain of $L^{1/2}$ coincides with the distributional fractional Sobolev space $\mathrm{S}^{\alpha,2}(\mathbb{R}^d)$ and that $\|L^{1/2} u\|_{\mathrm{L}^2} \simeq \|\nabla^\alpha u\|_{\mathrm{L}^2}$, with constants that are locally uniform in $\alpha$. As a consequence, we recover the Kato square root property for elliptic operators $-\operatorname{div}(\mu\nabla\cdot)$ in the limit $\alpha \to 1$, and we obtain continuity of the square root with respect to the parameter $\alpha$. The proof follows the second-order approach to the Kato problem. Its two new ingredients are nonlocal off-diagonal estimates of polynomial order for the resolvents of $L$, and a construction of $T(b)$-test functions that circumvents the nonlocality of $\nabla^\alpha$ by factorizing it through Riesz potentials.
\end{abstract}
 
\maketitle
 
\allowdisplaybreaks

\section{Introduction}

The main purpose of this article is to solve the Kato square root problem for fractional elliptic operators 
\begin{align}
\label{eq: operator L intro}
    L = -\div^\alpha(\mu \nabla^\alpha \cdot) \quad \text{on }\IR^d 
\end{align}
of order $2\alpha$, $\alpha \in (0,1)$, where $\nabla^\alpha$ and $\div^\alpha$ are the Riesz fractional gradient and divergence (see Definition~\ref{def: frac grad+div} below) and $\mu$ are bounded, measurable and uniformly elliptic complex coefficients. We show that the domain of $L^{\frac{1}{2}}$ coincides with the distributional fractional Sobolev space $\S^{\alpha,2}(\IR^d)$ consisting of all $\L^2(\IR^d)$-functions whose fractional gradient belongs to $\L^2(\IR^d;\IC^d)$, and that the homogeneous estimate $\|L^{\frac{1}{2}} u\|_{\L^2} \simeq \|\nabla^\alpha u\|_{\L^2}$ holds; see Theorem~\ref{thm: Kato property} below. 
 
\subsection*{The Kato square root problem}
 
At the beginning of the 1960s, Tosio Kato asked whether, given two Hilbert spaces $V \subseteq \H$ with $V$ dense in $\H$, the domain of the square root of the maximal accretive operator $L$ in $\H$ associated with a closed and sectorial sesquilinear form $\a :V \times V \to \IC$ always coincides with the form domain, that is, whether $\cD(L^{\frac{1}{2}}) = V$. Counterexamples by Lions~\cite{Lions} and McIntosh~\cite{McIntosh} showed that this fails in this generality. Lions already specialized the question to elliptic operators $L = - \div (\mu \nabla \cdot)$ with bounded measurable coefficients, and in this form it became known as the Kato square root problem. It was resolved in the affirmative on $\IR^d$ by Auscher, Hofmann, Lacey, McIntosh and Tchamitchian~\cite{AHLMcT} in 2002.
 
In the last two decades the square root property has been found to be deeply connected with boundary value problems~\cites{Auscher_Egert, AAAHK, Auscher_Axelsson_McIntosh2, AA1, AM, AusSta, HMiMo, AA}, operator-adapted function spaces~\cites{Auscher_Egert, Hofmann_Mayboroda_McIntosh, AusSta}, maximal regularity and quasilinear equations~\cites{AA1, ABHR, Pruss_Simonett}, and many further questions~\cites{Auscher, Lashi, HKMP, Bechtel, Egert}. Since the restriction from general sectorial forms to elliptic operators in divergence form is a drastic one, the question for which further classes of operators the square root property persists has been present throughout. For elliptic operators, the square root property was extended to rough boundary geometries and various boundary conditions~\cites{AKM_mixed, Egert_Haller_Tolksdorf, Bechtel_Egert_Haller, Bechtel-book, Bechtel_etal, Bechtel_Rosen} and to operators on submanifolds~\cite{Morris}. Using the theory of Muckenhoupt weights, it was extended to elliptic operators with degenerate coefficients~\cites{C-UR1, C-UR2, C-UMR}, and it was further established for Schrödinger operators~\cites{Bailey, Morris_Turner} and for parabolic operators~\cites{Auscher_Egert_Nystrom, Nystroem, Ataei_Egert_Nystrom}. As a key structural feature, all of these operators are built from classical differential operators, and hence local.
 
\subsection*{Nonlocality in the Kato problem}
 
Nonlocality has entered the Kato problem before. For example, Tolksdorf and the author~\cite{Haardt_Tolksdorf} established the square root property for the generalized Stokes operator, formally given by
\begin{align*}
    A u \coloneqq - \div(\mu \nabla u) + \nabla \phi, \quad \div(u) = 0 ,
\end{align*}
with bounded measurable coefficients $\mu$ satisfying a Gårding inequality. There, the nonlocality is caused by the pressure $\phi$ and the main difficulty it creates is the loss of local arguments, in particular of the classical exponential off-diagonal estimates for the resolvents~\cites{Tolksdorf_off-diagonal, Tolksdorf_caccioppoli}. A second, and rather different, form of nonlocality appears in the resolution of the parabolic Kato problem by Auscher, Egert and Nyström~\cite{Auscher_Egert_Nystrom}. The parabolic operator $\partial_t - \div_x(\mu \nabla_x \cdot)$ is itself local. However, if the operator is reformulated via a sesquilinear form or as a first-order system, nonlocal half-order derivatives in time enter, which create a major new challenge compared with the elliptic case. Fortunately, the nonlocality is confined to the time variable. This is a key observation made in \cite{Ataei_Egert_Nystrom}, where a second-order approach to the Kato property for weighted parabolic operators is established.\par

The operator $L=-\div^\alpha(\mu\nabla^\alpha\cdot)$ studied in this article, by contrast, is nonlocal in every direction. It is built from the Riesz fractional gradient
\begin{align}
\label{eq: fractional gradient}
    (\nabla^\alpha f)(x) = c_{d,\alpha} \int_{\IR^d}
    \frac{(y-x)\big(f(y)-f(x)\big)}{|x-y|^{d+1+\alpha}} \, \d y
    \qquad (x \in \IR^d)
\end{align}
and the fractional divergence
\begin{align*}
    (\div^\alpha F)(x) = c_{d,\alpha} \int_{\IR^d}
    \frac{(y-x) \cdot \big(F(y)-F(x)\big)}{|x-y|^{d+1+\alpha}} \, \d y
    \qquad (x \in \IR^d)
\end{align*}
of order $\alpha \in (0,1)$, where the normalization constant $c_{d,\alpha}$ is given in Definition~\ref{def: frac grad+div}. \v{S}ilhav\'y~\cite{Silhavy} singled out this pair as the unique fractional differential operators, uup to a multiplicative constant and subject to natural continuity requirements, that are invariant under translations and rotations and homogeneous of degree $\alpha$ under isotropic scaling. They are tied to the fractional Laplacian through $-\div^\alpha \nabla^\alpha = (-\Delta)^{\alpha}$, the operators $\nabla^\alpha$ and $-\div^\alpha$ are adjoint, and the classical derivatives are recovered as $\alpha \to 1$~\cites{Shieh_Spector1, Comi_Stefani1, Comi_StefaniI}. Moreover, they give rise to a distributional theory of fractional Sobolev spaces, which coincide with the Bessel potential spaces~\cites{Comi_Stefani1, Comi_StefaniII, Comi_Stefani2, Stefani-Thesis}, and to a theory of fractional partial differential equations~\cites{Shieh_Spector1, Shieh_Spector2, Schikorra_Shieh_Spector}. One motivation comes from nonlocal continuum mechanics, where these operators provide models, for example, of nonlocal hyperelasticity~\cite{Bellido_Cueto_MoraCorral_Piola} and, upon linearization, recover Eringen's model of nonlocal elasticity with singular kernels~\cite{Bellido_Cueto_MoraCorral_Eringen}. See also~\cites{Bellido_Cueto_MoraCorral, Cueto_Kreisbeck_Schoenberger, Kreisbeck_Schoenberger, Mazowiecka_Schikorra} for further developments.

\subsection*{Main results}
 
We impose the following assumption on the coefficients.
 
\begin{assumption}
\label{Ass: Coefficients}
    The coefficients $\mu \in \L^{\infty} (\IR^d ; \IC^{d\times d})$ satisfy for some $0 < \mu_{\bullet} \leq \mu^{\bullet}$ the inequalities
    \begin{align*}
        \Re \big(\mu(x) \xi \cdot \overline{\xi}\big) \geq \mu_\bullet |\xi|^2
        \quad \text{and} \quad
        |\mu(x)|_{\IC^{d\times d}} \leq \mu^{\bullet}
        \qquad (\xi \in \IC^d, \text{ a.e.\ } x \in \IR^d).
    \end{align*}
\end{assumption}
 
The operator $L$ in~\eqref{eq: operator L intro} is defined as the maximal accretive operator in $\L^2(\IR^d)$ associated with the sesquilinear form
\begin{align*}
    \a(u,v) = \int_{\IR^d} \mu \nabla^\alpha u \cdot \overline{\nabla^\alpha v} \, \d x
\end{align*}
on the distributional fractional Sobolev space $\S^{\alpha,2}(\IR^d)$; see Section~\ref{sec: second-order operator}. We emphasize that $L$ is in general not an integro-differential operator with a kernel comparable to $|x-y|^{-d-2\alpha}$, even for constant coefficients; see Remark~\ref{rem: kernel representation}.
Our main result reads as follows.
 
\begin{theorem}[The Kato square root property of $L$]
\label{thm: Kato property}
    Let $\alpha_0 \in (0,1)$ and let $\mu$ satisfy Assumption~\ref{Ass: Coefficients}. Then for all $\alpha\in[\alpha_0,1)$ we have $\cD(L^\frac{1}{2}) = \S^{\alpha,2}(\IR^d)$ and
    \begin{align*}
        \|L^{\frac{1}{2}} u\|_{\L^2} \simeq \|\nabla^\alpha u\|_{\L^2} \qquad (u \in \S^{\alpha,2}(\IR^d)),
    \end{align*}
    where the implicit constants depend only on $d$, $\alpha_0$, $\mu_\bullet$ and $\mu^\bullet$.
\end{theorem}

In particular, we obtain the Kato property for all $\alpha\in(0,1)$. Moreover, since $\S^{\alpha,2}(\IR^d)$ coincides with the fractional Sobolev--Slobodeckij space $\W^{\alpha,2}(\IR^d)$ (see Proposition~\ref{prop: homogeneous norms}), we obtain the following description.
 
\begin{corollary}
\label{cor: Kato2}
    Under the assumptions of Theorem~\ref{thm: Kato property},
    \begin{align*}
        \cD(L^\frac{1}{2}) =  \W^{\alpha,2}(\IR^d)
    \end{align*}
    and for all $u\in \W^{\alpha,2}(\IR^d)$ we have
    \begin{align*}
        \|L^\frac{1}{2} u\|_{\L^2} \simeq \|(-\Delta)^\frac{\alpha}{2} u\|_{\L^2}
        \simeq \alpha^\frac{1}{2}(1-\alpha)^\frac{1}{2}\|u\|_{\dot \W^{\alpha,2}}.
    \end{align*}
\end{corollary}

As consequences of the Kato property, we obtain the $\L^2$-boundedness of the Riesz transform $\nabla^\alpha L^{-\frac{1}{2}}$, as well as holomorphic and locally Lipschitz dependence of $L^{\frac{1}{2}}$ on the coefficients $\mu$ in the $\L^\infty$-topology; see Section~\ref{sec: application}.
Furthermore, a special feature of the Kato property in Theorem~\ref{thm: Kato property} is that the implicit constants are locally uniform with respect to the fractional parameter $\alpha$. In particular, our result fits into the asymptotic analysis of the fractional calculus as $\alpha\to 1$~\cites{Comi_StefaniI, Comi_StefaniII}, and we recover the classical Kato square root property for elliptic operators $L_1=-\div (\mu \nabla \cdot)$ in divergence form in the limit $\alpha\to 1$.

\begin{proposition}[The Kato property of $L_1$]
\label{prop: Kato for elliptic}
    Let $\mu$ satisfy Assumption~\ref{Ass: Coefficients}. Then we have $\cD(L_1^{\frac{1}{2}})=\W^{1,2}(\IR^d)$ and
    \begin{align*}
        \|L_1^{\frac{1}{2}}u\|_{\L^2}\simeq\|\nabla u\|_{\L^2}
        \qquad\big(u\in\W^{1,2}(\IR^d)\big),
    \end{align*}
    where the implicit constants depend only on $d$, $\mu_\bullet$ and $\mu^\bullet$.
\end{proposition}

Finally, we show continuity of the square root with respect to the parameter $\alpha$ for fixed coefficients $\mu$. To formulate such a statement, let us emphasize the fractional parameter $\alpha$ in the notation of $L$ by writing $L_{\alpha}$.

\begin{theorem}
\label{thm: continuity in alpha}
    Let $\mu$ satisfy Assumption~\ref{Ass: Coefficients}. Then for every $u\in\W^{1,2}(\IR^d)$ the map
    \begin{align*}
        (0,1]\to\L^2(\IR^d),\quad\alpha\mapsto L_\alpha^\frac{1}{2}u
    \end{align*}
    is continuous.
\end{theorem}

\subsection*{Strategy and main difficulties}

The main part of this paper is devoted to the Kato property in Theorem~\ref{thm: Kato property}. Here, we follow the second-order approach to the Kato problem for elliptic operators~\cite{AHLMcT}. By an abstract argument, Theorem~\ref{thm: Kato property} reduces to a square function estimate for $t^\alpha L(1+t^{2\alpha}L)^{-1}$, which is in turn reduced, by a principal part approximation, to a Carleson measure estimate, and the latter is proved by a $T(b)$ argument. In our setting, two steps of this scheme require new ideas.

\medskip

The first point concerns the off-diagonal estimates. Classical off-diagonal estimates for resolvents of elliptic operators possess an exponential decay rate and are classically derived through local arguments. To be more precise, the derivation relies on testing the resolvent equation with its solution $u$ multiplied by an exponential weight, or on an iteration of local Caccioppoli inequalities. However, due to the nonlocality of $\nabla^\alpha$ and the polynomial decay of its kernel, we can neither preserve exponential decay nor derive local Caccioppoli inequalities. A similar phenomenon also appeared in the treatment of the generalized Stokes operator~\cite{Haardt_Tolksdorf}. Nonetheless, it was observed there that the overall strategy for proving the Kato property remains stable under weaker nonlocal off-diagonal estimates with polynomial decay, provided the order is greater than the dimension $d$. For the Stokes operator, estimates of this type are derived via an iteration procedure involving nonlocal Caccioppoli inequalities~\cites{Haardt_Tolksdorf}. In our setting, we follow this approach by first establishing a novel nonlocal Caccioppoli inequality for the resolvent of $L$ in Theorem~\ref{thm: Caccioppoli} and then deriving the nonlocal off-diagonal estimates of polynomial order $\nu<d+\alpha_0$ by an iteration process in Theorem~\ref{thm: L2 off-diagonal estimates}.
Here, the main difficulty lies in controlling the nonlocal bilinear remainder $\nabla^\alpha_{\NL}$ coming from the product rule for $\nabla^\alpha$ (see Lemma~\ref{lem: Leibniz} and~\ref{lem: auxiliary lem}). Since $\alpha_0>0$, the decay order of our estimates can be chosen larger than $d$, which suffices for every step of the proof of Theorem~\ref{thm: Kato property}.\par
For comparison, in the first-order approach to the Kato property for the parabolic operator~\cite{Auscher_Egert_Nystrom}, the authors had to deal with similar issues due to the presence of a half-order time derivative. In particular, the off-diagonal estimates in~\cite{Auscher_Egert_Nystrom} are of restricted polynomial order and of anisotropic type due to the separation of space and time. More recently, Ataei, Egert and Nyström~\cite{Ataei_Egert_Nystrom} succeeded in providing a second-order proof of the Kato property for weighted parabolic operators, which makes use of the anisotropy so that off-diagonal estimates are only needed for operators that involve differentiation in space and no fractional operator. In this sense, the much more involved off-diagonal decay for nonlocal derivatives could be avoided. In our setting, however, the operator is nonlocal in every direction.

\medskip

The second difficulty concerns the Carleson measure argument, which relies on local $T(b)$-techniques. More precisely, for every cube $Q$, one has to construct a test function $b$ in the range of $\nabla^\alpha$ that points roughly in a given direction on $Q$ and satisfies $\|b\|_{\L^2} \lesssim |Q|^{\frac{1}{2}}$. In the elliptic setting, one chooses $b$ approximately equal to $\nabla \Phi$ with $\Phi$ being the product of an affine function, whose gradient points in the given direction, and a cut-off function that equals one on $Q$. Since $\nabla$ is a local operator, the cut-off does not affect the values and hence the direction of $\nabla \Phi$ on $Q$. For $\nabla^\alpha$, however, the values of $\nabla^\alpha \Phi$ on $Q$ depend on the values of $\Phi$ on all of $\IR^d$, so that the local construction no longer works.
Such a situation also appears in the parabolic setting~\cite{Auscher_Egert_Nystrom} due to the presence of the half-order time derivative. Luckily there, $\Phi$ has a product structure in space and time, and the time part is obtained by inverting the nonlocal half-order derivative on the Fourier side. This is possible because the time variable is one-dimensional, so that the Fourier symbol of the half-order derivative is scalar. In our setting, however, for $d \geq 2$ the Fourier symbol $m(\xi) = i\xi|\xi|^{\alpha-1}$ of $\nabla^\alpha$ is vector-valued and admits no right inverse. Thus, we instead exploit the factorization $\nabla^\alpha = \nabla I_{1-\alpha}$ through the Riesz potential, which formally gives $\nabla^\alpha (-\Delta)^{\frac{1-\alpha}{2}} = \nabla$. This reduces the construction to the local gradient $\nabla$, for which the elliptic construction is available again. We refer to Section~\ref{sec: Tb argument} for the details.

\medskip

Moreover, we emphasize that all steps in our procedure keep track of the dependence on the fractional parameter $\alpha$, which eventually leads to the local uniformity of the implicit constants in Theorem~\ref{thm: Kato property}. This will be the key to obtaining asymptotic results in the limit $\alpha\to 1$. Here, two effects have to be balanced: the normalization constant $c_{d,\alpha}$ in~\eqref{eq: fractional gradient} vanishes like $1-\alpha$, while integrals of the kernel $|x-y|^{-(d+\alpha)}$, such as in~\eqref{eq: kernel bound}, grow like $(1-\alpha)^{-1}$. In our estimates, both always appear together, so that they cancel out and the constants stay uniformly bounded with respect to $\alpha$. Moreover, we show that the fractional gradient $\nabla^\alpha$ converges to $\nabla$ and that the resolvents of $L=-\div^\alpha(\mu\nabla^\alpha\cdot)$ converge to the resolvents of $L_1=-\div(\mu\nabla\cdot)$. Together with the uniform square function estimate for $L$ obtained in the proof of Theorem~\ref{thm: Kato property}, this yields the corresponding square function estimate for $L_1$, and thus Proposition~\ref{prop: Kato for elliptic}. The same ingredients, combined with the locally uniform Kato bounds, give the continuity property in Theorem~\ref{thm: continuity in alpha}. Finally, the holomorphic and locally Lipschitz dependence on the coefficients in Proposition~\ref{prop: Holomorphic dependence} and Corollary~\ref{cor: Lipschitz estimate} follows as for elliptic operators~\cite{Auscher_Tchamitchian}, see also~\cite{Haardt_Tolksdorf}, since the argument only uses that $L$ is in divergence form and that the constants in Theorem~\ref{thm: Kato property} depend on $\mu$ only through $\mu_\bullet$ and $\mu^\bullet$.

\subsection*{Organization}

In Section~\ref{sec: fractional calculus} we collect the required properties of the Riesz fractional calculus that are used throughout the manuscript. In Section~\ref{sec: second-order operator} we introduce the fractional elliptic operator $L$ and prove a nonlocal Caccioppoli inequality and off-diagonal estimates for its resolvent. Section~\ref{sec: proof} contains the proofs of Theorem~\ref{thm: Kato property} and Corollary~\ref{cor: Kato2}. Finally, in Section~\ref{sec: application} we present applications of the Kato property, namely the $\L^2$-boundedness of the Riesz transform associated with $L$, the holomorphic and Lipschitz dependence on the coefficients, and the continuous dependence on the fractional parameter $\alpha$ as well as the limit $\alpha\to 1$.

\subsection*{Notation}

Throughout, $d \geq 1$ denotes the dimension and $\IN_0 = \IN \cup \{0\}$. We write $|\cdot|$ for the Euclidean norms on $\IC$ and $\IC^d$, and $a \cdot b = \sum_j a_j b_j$ for $a,b \in \IC^d$. The open ball with center $x$ and radius $r$ is denoted by $B=B(x,r)$, and for $\lambda > 0$ we write $\lambda B \coloneqq B(x,\lambda r)$ and $\lambda Q$ for the cube with the same center as $Q$ and $\lambda$ times its side-length $\ell(Q)$. For a ball $B$ and $k \in \IN_0$ we set $B_k \coloneqq 2^k B$ and define the dyadic annuli
\begin{align*}
    C_0(B) \coloneqq 2B, \qquad C_k(B) \coloneqq B_{k+1} \setminus \overline{B_{k}} \quad (k \geq 1),
\end{align*}
and analogously $C_k(Q)$ for cubes $Q$. For $\omega \in (0,\pi]$ let 
\begin{align*}
    \Sigma_{\omega} \coloneqq \{z \in \IC \setminus \{0\} : |\arg z| < \omega\}
\end{align*} 
and $\Sigma_{0} \coloneqq (0,\infty)$. We write $a \lesssim b$ if $a \leq C b$ for a constant $C > 0$ that depends only on parameters not at stake,  and $a \simeq b$ if $a\lesssim b$ and $b\lesssim a$ hold. In some situations it will be more handy to keep the notation $a \leq C b$ or $Ca \geq b$. In this case $C$ is generic.

For a measurable set $A$ we denote by $|A|$ its Lebesgue measure and by $\ind_A$ its indicator function. Moreover, if $A$ satisfies $0 < |A| < \infty$, then for $f \in \L^1_\loc(\IR^d)$ we write
\begin{align*}
    (f)_A \coloneqq \fint_A f \, \d x = \frac{1}{|A|} \int\limits_A f \, \d x.
\end{align*}
We denote by $\W^{1,2}(\IR^d)$ the classical Sobolev space. For $\beta \in (0 , 1)$ the fractional Sobolev--Slobodeckij space with differentiability $\beta$ will be denoted by $\W^{\beta,2}(\IR^d)$. We write its norm in integral form 
\begin{align*}
    \| f \|_{\W^{\beta,2}} \coloneqq \| f \|_{\L^2} + \| f \|_{\dot \W^{\beta,2}} \quad \text{with} \quad \| f \|_{\dot \W^{\beta,2}} \coloneqq \bigg(\iint_{\IR^d \times \IR^d} \frac{\lvert f(x) - f(y) \rvert^2}{\lvert x - y \rvert^{d + 2 \beta}} \, \d x \, \d y \bigg)^{\frac{1}{2}}.
\end{align*}
For $f \in \W^{\beta,2} (\IR^d)$ we will need to calculate a counterpart of the homogeneous $\dot \W^{\beta,2}$-seminorm on a measurable set $A \subseteq \IR^d$ which we write as
\begin{align*}
 \| f \|_{\dot \W^{\beta,2} (A)} \coloneqq \bigg(\iint_{A \times A} \frac{\lvert f(x) - f(y) \rvert^2}{\lvert x - y \rvert^{d + 2 \beta}} \, \d x \, \d y \bigg)^{\frac{1}{2}}.
\end{align*}

\section{A glimpse into fractional calculus}
\label{sec: fractional calculus}

In this section we introduce the fractional calculus of Šilhavý~\cite{Silhavy} and collect the properties of its operators that are needed later on, see also~\cites{Shieh_Spector1, Comi_Stefani1, Comi_StefaniII}.
 
We begin with the integral representations of the fractional gradient and the fractional divergence.

\begin{definition}
\label{def: frac grad+div}
    Let $\alpha \in (0,1)$ and set
    \begin{align*}
        c_{d,\alpha} \coloneqq 2^\alpha \pi^{-\frac{d}{2}}
        \frac{\Gamma \big(\frac{d+\alpha+1}{2} \big)}{\Gamma \big(\frac{1-\alpha}{2} \big)}.
    \end{align*}
    The \textit{$\alpha$-fractional gradient} of $f \in \C_c^\infty(\IR^d)$ and the
    \textit{$\alpha$-fractional divergence} of $F \in \C_c^\infty(\IR^d;\IC^d)$ are
    \begin{align*}
        (\nabla^\alpha f)(x) \coloneqq c_{d,\alpha} \int\limits_{\IR^d}
        \frac{(z-x)\big(f(z)-f(x)\big)}{|x-z|^{d+\alpha+1}} \, \d z
    \end{align*}
    and 
    \begin{align*}
        (\div^\alpha F)(x) \coloneqq c_{d,\alpha} \int\limits_{\IR^d}
        \frac{(z-x) \cdot \big(F(z)-F(x)\big)}{|x-z|^{d+\alpha+1}} \, \d z
    \end{align*}
    for $x \in \IR^d$. 
\end{definition}

\begin{remark}
    The constant $c_{d,\alpha}$ is chosen such that $\nabla^\alpha f \to \nabla f$ and $\div^\alpha F \to \div F$ as $\alpha \to 1$, see \cite[Thm.~4.11]{Comi_StefaniI}, and satisfies
    \begin{align}
    \label{eq: normalization constant}
        \frac{c_{d,\alpha}}{1-\alpha} \leq C(d)
    \end{align}
    for all $\alpha \in (0,1)$ and some constant $C(d)>0$ depending only on the dimension $d$, see \cite[Lem.~4.1]{Comi_StefaniI}.
\end{remark}

Define the Fourier transform
\begin{align*}
    \cF(f)(\xi) \coloneqq \frac{1}{(2\pi)^{\frac{d}{2}}} \int\limits_{\IR^d}
    f(x) e^{-i x \cdot \xi} \, \d x.
\end{align*}
The following spectral representation is due to Shieh and Spector~\cite[Thm.~1.4]{Shieh_Spector1}.

\begin{proposition}
\label{prop: fourier symbols}
    Let $\alpha \in (0,1)$. For all $f \in \C_c^\infty(\IR^d)$ and all
    $F \in \C_c^\infty(\IR^d;\IC^d)$ we have
    \begin{align*}
        \cF(\nabla^\alpha f)(\xi) = \frac{i \xi}{|\xi|^{1-\alpha}} \cF(f)(\xi)
        \quad \text{and} \quad
        \cF(\div^\alpha F)(\xi) = \frac{i \xi^\top}{|\xi|^{1-\alpha}} \cF(F)(\xi)
        \qquad (\xi \in \IR^d) .
    \end{align*}
    In particular, we have
    \begin{align*}
        -\div^\alpha \nabla^\alpha = (-\Delta)^{\alpha},
    \end{align*}
    where $(-\Delta)^{\alpha} = \cF^{-1}(|\xi|^{2\alpha}\cF)$ denotes the fractional Laplacian.
\end{proposition}

\subsection{Distributional fractional Sobolev spaces}

As for the classical operators, $\nabla^\alpha$ and $-\div^\alpha$ are adjoint to each other, see~\cite[Lem.~2.5]{Comi_Stefani1}. More precisely, for $\alpha \in (0,1)$,  $f \in \C_c^\infty(\IR^d)$ and $F\in \C_c^\infty(\IR^d;\IC^d)$ we have
\begin{align*}
    \int\limits_{\IR^d} \nabla^\alpha f \cdot F \, \d x =-\int\limits_{\IR^d} f \, \div^\alpha F \, \d x.
\end{align*}
This duality relation is the starting point for a distributional fractional calculus imitating the classical Sobolev spaces.

\begin{definition}
\label{def: S space}
    Let $\alpha \in (0,1)$. Given $f \in \L^2(\IR^d)$, we say that $g \in \Lloc^1(\IR^d;\IC^d)$ is a \textit{weak $\alpha$-gradient} of $f$, written $g = \nabla^\alpha f$, if
    \begin{align*}
        \int\limits_{\IR^d} g \cdot F \, \d x
        = -\int\limits_{\IR^d} f \, \div^\alpha F \, \d x
        \qquad (F \in \C_c^\infty(\IR^d;\IC^d)).
    \end{align*}
    We set
    \begin{align*}
        \S^{\alpha,2}(\IR^d) \coloneqq
        \big\{ f \in \L^2(\IR^d) : \nabla^\alpha f \in \L^2(\IR^d;\IC^d) \big\}\quad \text{with} \quad
        \|f\|_{\S^{\alpha,2}} \coloneqq \|f\|_{\L^2} + \|\nabla^\alpha f\|_{\L^2},
    \end{align*}
    and we write $\S^{-\alpha,2}(\IR^d) \coloneqq \big(\S^{\alpha,2}(\IR^d)\big)'$ for the space of bounded antilinear functionals on $\S^{\alpha,2}(\IR^d)$.
\end{definition}

The space $\S^{\alpha,2}(\IR^d)$ is the natural domain of the multipliers of Proposition~\ref{prop: fourier symbols}, and can be identified with the fractional Sobolev space $\W^{\alpha,2}(\IR^d)$.

\begin{proposition}
\label{prop: homogeneous norms}
    Let $\alpha \in (0,1)$. Then $\S^{\alpha,2}(\IR^d) 
        = \W^{\alpha,2}(\IR^d)$
    with equivalent norms, $\C_c^\infty(\IR^d)$ is dense in $\S^{\alpha,2}(\IR^d)$, and for all $u\in\S^{\alpha,2}(\IR^d)$ we have
    \begin{align}
    \label{eq: homogeneous norms}
        \|\nabla^\alpha u\|_{\L^2}^2 = \|(-\Delta)^{\frac{\alpha}{2}} u\|_{\L^2}^2 \simeq \alpha(1-\alpha) \|u\|^2_{\dot{\W}^{\alpha,2}}
    \end{align}
    with implicit constants depending only on $d$. 
\end{proposition}

\begin{proof}
    First, density of $\C_c^\infty(\IR^d)$ in $\W^{\alpha,2}(\IR^d)$ and $\S^{\alpha,2}(\IR^d)$ was proven in \cite[Thm.~2.4]{Hitchhiker} and \cite[Thm.~25]{Comi_StefaniII}, respectively. Now, let $u \in \C_c^\infty(\IR^d)$. Since $|i \xi |\xi|^{\alpha-1}| = |\xi|^\alpha$, Proposition~\ref{prop: fourier symbols} and Plancherel's theorem give the first equality in~\eqref{eq: homogeneous norms}. For the second equivalence, by \cite[Prop.~3.6]{Hitchhiker}, we have
    \begin{align*}
        \|(-\Delta)^{\frac{\alpha}{2}} u\|^2_{\L^2} = \frac{C(d,\alpha) }{2}\|u\|^2_{\dot{\W}^{\alpha,2}},
    \end{align*}
    where the constant $C(d,\alpha)>0$ has the closed form
    \begin{align*}
        C(d,\alpha) = \bigg(\int\limits_{\IR^d} \frac{1-\cos(\xi_1)}{|\xi|^{d+2\alpha}}\,\d \xi\bigg)^{-1} = \frac{\alpha 4^{\alpha} \Gamma(\frac{d}{2}+\alpha)}{\pi^\frac{d}{2}\Gamma(1-\alpha)},
    \end{align*}
    see also \cite[Sec.~1.10]{Ros-Oton}. Using the identity $(1-\alpha)\Gamma(1-\alpha) = \Gamma(2-\alpha)$, we get
    \begin{align*}
        C(d,\alpha)= \alpha (1-\alpha) \frac{4^{\alpha} \Gamma(\frac{d}{2}+\alpha)}{\pi^\frac{d}{2}\Gamma(2-\alpha)}.
    \end{align*}
    Since the last factor on the right-hand side is continuous and positive in $\alpha\in[0,1]$, there exist implicit constants depending only on $d$ such that $C(d,\alpha)\simeq \alpha (1-\alpha)$.
    Finally, due to density of $\C_c^\infty(\IR^d)$ in both fractional spaces, it follows that $\S^{\alpha,2}(\IR^d) = \W^{\alpha,2}(\IR^d)$ with equivalent norms and that \eqref{eq: homogeneous norms} extends to $u\in \S^{\alpha,2}(\IR^d)$.
\end{proof}

In particular, we can extend the spectral representation of $\nabla^\alpha$ from Proposition~\ref{prop: fourier symbols} to $\S^{\alpha,2}(\IR^d)$ via density. We will use this fact freely, in particular in the next statement to obtain the following embedding chain; see also \cite[Prop.~33]{Comi_StefaniII}.

\begin{lemma}
\label{lem: embedding}
    Let $0<\beta<\alpha <1$. Then we have the continuous embeddings
    \begin{align*}
        \W^{1,2}(\IR^d) \hookrightarrow \S^{\alpha,2}(\IR^d) \hookrightarrow \S^{\beta,2}(\IR^d).
    \end{align*}
\end{lemma}

\begin{proof}
    Since $1+|\xi|^\beta \leq 2+|\xi|^\alpha \leq 3+|\xi|$, we obtain with Proposition~\ref{prop: fourier symbols} and Plancherel's theorem
    \begin{align*}
        \|u\|_{\S^{\beta,2}} \lesssim \|u\|_{\S^{\alpha,2}} \lesssim \|u\|_{\W^{1,2}}
    \end{align*}
    for $u\in \C_c^\infty(\IR^d)$. Since $\C_c^\infty(\IR^d)$ is a dense subspace in all of the three spaces, we conclude the proof.
\end{proof}

An alternative representation of the fractional gradient can be given in terms of the classical gradient and the Riesz potential. More precisely, for $\beta \in (0,1)$, define the Riesz potential $I_\beta$ on $\C_c^\infty(\IR^d)$ by
\begin{align}
\label{eq: def Riesz potential}
    (I_\beta f)(x) \coloneqq \frac{c_{d,1-\beta}}{d-\beta}
    \int\limits_{\IR^d} \frac{f(y)}{|x-y|^{d-\beta}} \, \d y
    \qquad (x \in \IR^d).
\end{align}
Then, we have
\begin{align}
\label{eq: frac op via Riesz potential}
    \nabla^\alpha f= I_{1-\alpha}\nabla f= \nabla I_{1-\alpha}f,
\end{align}
for $f\in \C_c^\infty(\IR^d)$, see \cite[Prop.~2.2]{Comi_Stefani1}.
Since $(-\Delta)^{\frac{1-\alpha}{2}}$ is the inverse of $I_{1-\alpha}$, see~\cite[Sec.~2.3]{Comi_Stefani1} and~\cite[Sec.~V.1]{Stein}, composing the identity above with $(-\Delta)^{\frac{1-\alpha}{2}}$ restores (formally) the local gradient
\begin{align*}
    \nabla^\alpha (-\Delta)^{\frac{1-\alpha}{2}}f
    = \nabla I_{1-\alpha} (-\Delta)^{\frac{1-\alpha}{2}}f = \nabla f.
\end{align*}
This relation between the classical and fractional gradient in the context of Sobolev functions was shown in~\cite[Prop.~3.1 and~3.3]{Kreisbeck_Schoenberger} and is recorded in the following proposition.

\begin{proposition}
\label{prop: local correspondence}
    Let $\alpha \in (0,1)$. For every $v \in \W^{1,2}(\IR^d)$ the function
    \begin{align*}
        u \coloneqq (-\Delta)^{\frac{1-\alpha}{2}} v
    \end{align*}
    belongs to $\S^{\alpha,2}(\IR^d)$ and satisfies $\nabla^\alpha u = \nabla v$.
\end{proposition}

We also need a fractional divergence that is defined on all of $\L^2(\IR^d;\IC^d)$.

\begin{definition}
\label{def: very weak div}
    Let $\alpha \in (0,1)$. The \textit{very weak fractional divergence}
    $\div^\alpha :\L^2(\IR^d;\IC^d) \to \S^{-\alpha,2}(\IR^d)$ is defined by
    \begin{align*}
        \langle \div^\alpha F , g \rangle
        \coloneqq -\int\limits_{\IR^d} F \cdot \overline{\nabla^\alpha g} \, \d x
        \qquad \big(g \in \S^{\alpha,2}(\IR^d)\big) .
    \end{align*}
\end{definition}

\subsection{The Leibniz rule}

In contrast to the classical case, the product rule for $\nabla^\alpha$ carries a genuinely nonlocal, bilinear remainder. More precisely, for $f,g\in  \C_c^\infty(\IR^d)$ we define
\begin{align*}
    \nabla^\alpha_{\NL}(f,g)(x)
    \coloneqq c_{d,\alpha} \int_{\IR^d}
    \frac{(y-x)(f(y)-f(x))(g(y)-g(x))}
    {|x-y|^{d+\alpha+1}} \, \d y.
\end{align*}
It is this remainder that has to be controlled in the off-diagonal estimates of Section~\ref{sec: caccioppoli and off-diag est}.
The following Leibniz rule for the fractional gradient can be found in~\cite[Lem.~2.11]{Kreisbeck_Schoenberger}; see also \cite[Thm.~4.1]{Comi_Stefani2}.

\begin{lemma}[Leibniz rule]
\label{lem: Leibniz}
    Let $\alpha \in (0,1)$. For $\varphi \in \C_c^\infty(\IR^d)$ and $u\in \S^{\alpha,2}(\IR^d)$ we have $u\varphi \in \S^{\alpha,2}(\IR^d)$ and
    \begin{align*}
        \nabla^\alpha (u\varphi)
        = \varphi\nabla^\alpha u + u \nabla^\alpha \varphi
        + \nabla^\alpha_{\NL}(u,\varphi),
    \end{align*}
    where all terms on the right-hand side belong to $\L^2(\IR^d;\IC^d)$.
\end{lemma}

We end this section by providing estimates for $\nabla^\alpha$ and the remainder $\nabla^\alpha_\NL$ that are needed in Sections~\ref{sec: second-order operator} and~\ref{sec: proof}.

\begin{lemma}
\label{lem: auxiliary lem}
    Let $\alpha_0\in(0,1)$, $B=B(x_0,r)$ be a ball, $\rho\in(0,r]$ and let $\varphi\in\C_c^\infty(\IR^d)$ satisfy $0\leq\varphi\leq 1$, $\supp\varphi\subset2B$ and $\rho \|\nabla\varphi\|_{\L^\infty}\leq M$ for some $M>0$. Then there exists a constant $C>0$ depending only on $d,M$ and $\alpha_0$ such that for all $\alpha\in[\alpha_0,1)$ the following holds.
    \begin{enumerate}[label=(\roman*)]
        \item\label{item: est1} For every $k\in\IN_0$ we have
        \begin{align*}
            \|\nabla^\alpha\varphi\|_{\L^\infty(C_k(B))}
            \leq C\,\frac{2^{-(d+\alpha)k}}{\rho^{\alpha}} .
        \end{align*}
        \item\label{item: est2} For $u\in\L^2(\IR^d)$ we have
        \begin{align*}
            \|\nabla^\alpha_{\NL}(u,\varphi)\|_{\L^2(2B)}
            \leq \frac{C}{\rho^\alpha}\sum_{m=0}^\infty2^{-(\frac d2+\alpha)m}\|u\|_{\L^2(C_m(B))}
        \end{align*}
        and, for every $k\geq 1$
        \begin{align*}
            \|\nabla^\alpha_{\NL}(u,\varphi)\|_{\L^2(C_k(B))}
            \leq\frac{C}{\rho^\alpha}\Big(2^{-(d+\alpha)k}\|u\|_{\L^2(C_k(B))}
            +2^{-(\frac d2+\alpha)k}\|u\|_{\L^2(2B)}\Big).
        \end{align*}
    \end{enumerate}
\end{lemma}

\begin{proof}
    Since $0\leq\varphi\leq1$ and $\rho\|\nabla\varphi\|_{\L^\infty}\lesssim1$, we have $|\varphi(x)-\varphi(y)|\lesssim\min\{|x-y|/\rho,1\}$. Set $\kappa(z)\coloneqq\min\{|z|/\rho,1\}\,|z|^{-d-\alpha}$ and notice that
    \begin{align}
    \label{eq: kernel bound}
        \|\kappa\|_{\L^1} = \frac{1}{\rho}\int\limits_{B(0,\rho)}\frac{1}{|z|^{d+\alpha-1}}\,\d z +\int\limits_{\IR^d\setminus B(0,\rho)}\frac{1}{|z|^{d+\alpha}}\,\d z \lesssim 
        \frac{1}{(1-\alpha)\rho^\alpha}.
    \end{align}
    We divide the proof into estimating the fractional gradient and the bilinear term.\\
    
    \noindent \textbf{Step 1: Show \ref{item: est1}}.
    For every $x\in\IR^d$ we calculate
    \begin{align*}
        |(\nabla^\alpha\varphi)(x)| &\lesssim c_{d,\alpha}\int\limits_{\IR^d} |\kappa(x-y)|\,\d y  = c_{d,\alpha}\|\kappa\|_{\L^1} \lesssim \frac{ 1}{\rho^\alpha},
    \end{align*}
    where we used \eqref{eq: kernel bound} and the uniform estimate \eqref{eq: normalization constant}. This proves \ref{item: est1} for $k=0,1$.
    If  $x\in C_k(B)$, for $k\geq 2$, we notice that $\varphi(x) = 0$. Thus,
    \begin{align*}
        |(\nabla^\alpha\varphi) (x)| \leq c_{d,\alpha}\int\limits_{\IR^d} \frac{|\varphi(y)- \varphi(x)|}{|x-y|^{d+\alpha}}\,\d y =  c_{d,\alpha}\int\limits_{2B} \frac{|\varphi(y)|}{|x-y|^{d+\alpha}}\,\d y.
    \end{align*}
    Moreover, for $x\in C_k(B)$ and $y\in 2B$ we have $|x-y|\gtrsim 2^{k}r$ so that with \eqref{eq: normalization constant} and $\rho\leq r$ it follows
    \begin{align*}
        |(\nabla^\alpha\varphi) (x)| \lesssim \frac{c_{d,\alpha}}{(2^kr)^{d+\alpha}}\int\limits_{2B} |\varphi(y)|\,\d y\lesssim \frac{2^{-(d+\alpha)k}}{\rho^{\alpha}}.
    \end{align*}

    \noindent \textbf{Step 2: Show \ref{item: est2}.}
    Notice that for all $x\in\IR^d$ we have
    \begin{align*}
        |\nabla^\alpha_{\NL}(u,\varphi)(x)|
        &\leq c_{d,\alpha}\int\limits_{\IR^d}\frac{|u(x)-u(y)| |\varphi(x)-\varphi(y)|}{|x-y|^{d+\alpha}}\,\d y \\
        &\leq c_{d,\alpha}\bigg(|u(x)|\int\limits_{\IR^d}\frac{ |\varphi(x)-\varphi(y)|}{|x-y|^{d+\alpha}}\,\d y +\int\limits_{\IR^d}\frac{|u(y)| |\varphi(x)-\varphi(y)|}{|x-y|^{d+\alpha}}\,\d y\bigg).
    \end{align*}
    Thus,  by Step~1, for every $k\in\IN_0$ we have the bound
    \begin{align}
    \label{eq: est on NL}
        \|\nabla^\alpha_{\NL}(u,\varphi)\|_{\L^2(C_k(B))}&\lesssim \frac{2^{-(d+\alpha)k}\|u\|_{\L^2(C_k(B))}}{\rho^\alpha} + c_{d,\alpha}\bigg(\int\limits_{C_k(B)}\bigg(\int\limits_{\IR^d}\frac{|u(y)| |\varphi(x)-\varphi(y)|}{|x-y|^{d+\alpha}}\,\d y\bigg)^2\,\d x\bigg)^\frac{1}{2}\notag\\
        &\lesssim \frac{2^{-(d+\alpha)k}\|u\|_{\L^2(C_k(B))}}{\rho^\alpha} + c_{d,\alpha}\sum\limits_{m=0}^\infty \I(k,m),
    \end{align}
    where we set
    \begin{align*}
        \I(k,m)\coloneqq \bigg(\int\limits_{C_k(B)}\bigg(\int\limits_{C_m(B)}\frac{|u(y)| |\varphi(x)-\varphi(y)|}{|x-y|^{d+\alpha}}\,\d y\bigg)^2\,\d x\bigg)^\frac{1}{2} \qquad (k,m\in\IN_0).
    \end{align*}
    We further investigate $\I(k,m)$ by distinguishing different cases.\\

    \noindent \textbf{Case 1: $k=0$.}
    Assume $0\leq m\leq 2$. We use Young's convolution inequality and \eqref{eq: kernel bound} to obtain
    \begin{align*}
        \I(0,m)\lesssim \|\kappa \ast (\ind_{C_m(B)}u)\|_{\L^2(2B)}\leq\|\kappa\|_{\L^1}\|u\|_{\L^2(C_m(B))} \lesssim\frac{\|u\|_{\L^2(C_m(B))}}{(1-\alpha)\rho^\alpha}.
    \end{align*}
    In the case $m>2$, we notice that $\varphi(y) = 0$ and $|x-y|\gtrsim 2^{m}r$ for $x\in 2B$ and $y\in C_m(B)$. Hence, we estimate with Hölder's inequality
    \begin{align*}
        \I(0,m) \lesssim  \bigg(\int\limits_{2B}\bigg(\int\limits_{C_m(B)}\frac{|u(y)| |\varphi(x)|}{(2^mr)^{d+\alpha}}\,\d y\bigg)^2\,\d x\bigg)^\frac{1}{2}
        &\leq  \frac{|2B|^\frac{1}{2}}{(2^mr)^{d+\alpha}}\|u\|_{\L^1(C_m(B))} \\
        &\leq\frac{|2B|^\frac{1}{2}|C_m(B)|^\frac{1}{2}}{(2^mr)^{d+\alpha}}\|u\|_{\L^2(C_m(B))} \\
        &\lesssim \frac{2^{-(\frac{d}{2}+\alpha)m}}{\rho^\alpha}\|u\|_{\L^2(C_m(B))},
    \end{align*}
    where we used $\rho\leq r$ in the last step. Combining both estimates in \eqref{eq: est on NL}, and using the uniform bound \eqref{eq: normalization constant}, we obtain
    \begin{align*}
        \|\nabla^\alpha_{\NL}(u,\varphi)\|_{\L^2(2B)}
        \lesssim\sum\limits_{m=0}^{\infty} \frac{2^{-(\frac{d}{2}+\alpha)m}}{\rho^\alpha}\|u\|_{\L^2(C_m(B))}.
    \end{align*}

    \noindent \textbf{Case 2: $k\geq 1$.}
    If $m\geq1$, we get $\I(k,m)=0$ because $\varphi=0$ on $C_m(B)$ and $C_k(B)$. Thus, assume $m=0$. If $k=1$, we use Young's convolution inequality and \eqref{eq: kernel bound} to arrive at
    \begin{align*}
        \I(1,0)\lesssim \|\kappa\ast (\ind_{2B}u)\|_{\L^2(C_1(B))}\lesssim\|\kappa\|_{\L^1}\|u\|_{\L^2(2B)}\lesssim\frac{\|u\|_{\L^2(2B)}}{(1-\alpha)\rho^\alpha}.
    \end{align*}
    If $k\geq2$, then $\varphi(x)=0$ and $|x-y|\gtrsim 2^{k}r$ for $x\in C_k(B)$ and $y\in 2B$, so that
    \begin{align*}
        \I(k,0) &\lesssim  \bigg(\int\limits_{C_k(B)}\bigg(\int\limits_{2B}\frac{|u(y)| |\varphi(y)|}{(2^kr)^{d+\alpha}}\,\d y\bigg)^2\,\d x\bigg)^\frac{1}{2}\\
        &\leq  \frac{|C_k(B)|^\frac{1}{2}}{(2^kr)^{d+\alpha}}\|u\|_{\L^1(2B)} \\
        &\leq\frac{|C_k(B)|^\frac{1}{2}|2B|^\frac{1}{2}}{(2^kr)^{d+\alpha}}\|u\|_{\L^2(2B)} \\
        &\lesssim \frac{2^{-(\frac{d}{2}+\alpha)k}}{\rho^\alpha}\|u\|_{\L^2(2B)},
    \end{align*}
    where we used Hölder's inequality in the second last step and $\rho\leq r$ in the last step. Combining both cases in \eqref{eq: est on NL}, and using the uniform bound \eqref{eq: normalization constant}, we obtain the estimate
    \begin{align*}
        \|\nabla^\alpha_{\NL}(u,\varphi)\|_{\L^2(C_k(B))}
        \lesssim \frac{2^{-(d+\alpha)k}\|u\|_{\L^2(C_k(B))}}{\rho^\alpha} + \frac{2^{-(\frac{d}{2}+\alpha)k}}{\rho^\alpha}\|u\|_{\L^2(2B)}.
    \end{align*}
    This yields \ref{item: est2} for every $u\in\L^2(\IR^d)$ in the case $k\geq 1$.
\end{proof}

\begin{remark}
\label{rem: cubes instead of balls}
    Lemma~\ref{lem: auxiliary lem} can be stated and proved for cubes $Q$ instead of balls $B$ in a straightforward way, where $r$ is replaced by $\ell(Q)$.
\end{remark}

\section{Fractional elliptic operator}
\label{sec: second-order operator}

In this section, we define fractional elliptic operators in divergence form and collect important properties of them. Throughout the section, let $\alpha\in (0,1)$ and $\mu\in \L^\infty(\IR^d;\IC^{d\times d})$ satisfy Assumption~\ref{Ass: Coefficients}. Moreover, we set
\begin{align}
\label{def: omega_0}
    \omega_0\coloneqq \arctan\Big(\frac{\mu^\bullet}{\mu_\bullet}\Big) \in \Big(0,\frac{\pi}{2}\Big).
\end{align}

\begin{definition}
\label{def: ses-form and operator}
    We define the sesquilinear form
    \begin{align*}
        \a : \S^{\alpha,2}(\IR^d) \times \S^{\alpha,2}(\IR^d) \to \IC, \quad
        \a(u,v) \coloneqq \int\limits_{\IR^d}\mu\nabla^\alpha u \cdot
        \overline{\nabla^\alpha v} \, \d x
    \end{align*}
    and the associated \textit{weak fractional operator} $\cL:\S^{\alpha,2}(\IR^d)\to \S^{-\alpha,2}(\IR^d)$ by
    \begin{align*}
        \langle \cL u,v\rangle \coloneqq \a(u,v) \qquad (u,v\in\S^{\alpha,2}(\IR^d)).
    \end{align*}
    Subsequently, we define the fractional operator $L$ as the maximal restriction of $\cL$ to $\L^2(\IR^d)$, that is
    \begin{align*}
        \cD(L) = \{u\in \S^{\alpha,2}(\IR^d): \cL u\in \L^2(\IR^d)\} \quad \text{with} \quad Lu = \cL u. 
    \end{align*}
\end{definition}

\begin{remark}
\label{rem: kernel representation}
    The operator $L$ is in general not an integro-differential operator with a kernel comparable to $|x-y|^{-d-2\alpha}$, even for constant coefficients. Indeed, if $\mu$ is real, symmetric
    and constant, then there exists a unique symmetric kernel $k_\mu\in\Lloc^1(\IR^{2d}\setminus\{x=y\})$
    such that
    \begin{align*}
        \a(u,v) = \iint\limits_{\IR^d\times\IR^d} \big(u(x)-u(y)\big)\overline{\big(v(x)-v(y)\big)}\,k_\mu(x,y)\,\d x\,\d y
        \qquad \big(u,v\in\C_c^\infty(\IR^d)\big),
    \end{align*}
    and $k_\mu(x,y)|x-y|^{d+2\alpha}$ is bounded from below by a positive constant if and only if the smallest eigenvalue of $\mu$, denoted by $\lambda_{\min}(\mu)$, satisfies
    $(d+2\alpha)\lambda_{\min}(\mu)>\operatorname{tr}(\mu)$; see~\cite[Thm.~1.2]{Dipierro_Sportelli_Valdinoci}. We also refer to~\cites{Carbotti,Lo_Rodrigues} for related questions.
\end{remark}

For later reasoning it will be important to deduce a representation of $\cL$ as a composition of the fractional gradient $\nabla^\alpha$, the coefficients $\mu$ and the very weak fractional divergence $\div^\alpha$ from Definition~\ref{def: very weak div}.

\begin{lemma}
\label{lem: identity for cL}
    For all $u\in \S^{\alpha,2}(\IR^d)$, we have $\cL u = -\div^\alpha(\mu\nabla^\alpha u)$. 
\end{lemma}

Under Assumption~\ref{Ass: Coefficients}, the sesquilinear form $\a$ is bounded and satisfies
\begin{align*}
    \a(u,u) \in \overline{\Sigma_{\omega_0}} \quad \text{and} \quad \Re \a (u,u)\geq \mu_\bullet \|\nabla^\alpha u\|_{\L^2}^2
\end{align*}
for all $u\in\S^{\alpha,2}(\IR^d)$. Hence, classical form theory yields that $L$ is maximal accretive and sectorial of angle $\omega_0$ on $\L^2(\IR^d)$. In particular, $L$ possesses a unique maximal accretive square root $L^\frac{1}{2}$, which is the operator appearing in Theorem~\ref{thm: Kato property}. We refer to~\cites{Kato, ISem27} for more details. Finally, we have the following uniform resolvent bounds.

\begin{proposition}
\label{prop: L2 resolvent bounds}
    We have $\sigma(L) \cup \sigma(\cL) \subset
    \overline{\Sigma_{\omega_0}}$ and for all $\theta\in [0,\pi-\omega_0)$ there exists $C>0$ such that for all $\lambda \in \Sigma_\theta$ and all $f\in \L^2(\IR^d)$ we have
    \begin{align*}
        |\lambda|\|(\lambda + L)^{-1}f\|_{\L^2}
        + |\lambda|^\frac{1}{2}\|\nabla^\alpha (\lambda + L)^{-1}f\|_{\L^2}
        \leq C\|f\|_{\L^2}.
    \end{align*}
    Moreover, there exists $C>0$ such that for all $\lambda \in \Sigma_\theta$ and all $F\in \L^2(\IR^d;\IC^{d})$ we have
    \begin{align*}
        |\lambda|^\frac{1}{2}\|(\lambda + \cL)^{-1}\div^\alpha F\|_{\L^2}
        + \|\nabla^\alpha (\lambda + \cL)^{-1}\div^\alpha F\|_{\L^2}
        \leq C\|F\|_{\L^2}.
    \end{align*}
    The constant $C$ depends only on $d$, $\theta$ and $\mu_\bullet, \mu^{\bullet}$.
\end{proposition}

The proof is a standard application of the Lax--Milgram lemma, and will therefore be omitted.

Finally, for $x\in\IR^d$ let $\mu^*(x)$ denote the adjoint of $\mu(x)$. Then, for all $u,v\in\S^{\alpha,2}(\IR^d)$ we have
\begin{align*}
    \int\limits_{\IR^d} \mu^*\nabla^\alpha u\cdot\overline{\nabla^\alpha v} \,\d x
    = \overline{\int\limits_{\IR^d} \mu\nabla^\alpha v\cdot\overline{\nabla^\alpha u }\, \d x},
\end{align*}
so that $\mu^*$ again satisfies Assumption~\ref{Ass: Coefficients} with the same constants $\mu_\bullet, \mu^\bullet$. Thus, introducing the adjoint form
\begin{align*}
    \a^*(u,v) \coloneqq \overline{\a(v,u)}
    = \int\limits_{\IR^d} \mu^*\nabla^\alpha u\cdot\overline{\nabla^\alpha v}\, \d x
    \qquad (u,v\in\S^{\alpha,2}(\IR^d)),
\end{align*}
we denote by $\cL^*$ and $L^*$ the (weak) fractional operators associated with $\a^*$. These are the adjoints of $\cL$ and $L$, have the same sectoriality angle $\omega_0$, and every result of this paper applies verbatim to them. We use this in duality arguments without further mention.

\subsection{Caccioppoli inequality and off-diagonal estimates}
\label{sec: caccioppoli and off-diag est}

In this subsection, we follow the approach for the generalized Stokes operator from \cites{Haardt_Tolksdorf, Tolksdorf_caccioppoli} to first derive a nonlocal Caccioppoli inequality (Theorem~\ref{thm: Caccioppoli}) and then deduce off-diagonal estimates by an iteration scheme (Theorem~\ref{thm: L2 off-diagonal estimates}). We start with the following lemma.

\begin{lemma}
\label{lem: iteration setup}
    Let $\alpha_0\in(0,1)$. For all $\theta \in [0,\pi-\omega_0)$ and all $\nu \in (0, d+\alpha_0)$ there exists a constant $C >0$ such that for all balls $B= B(x_0,r)$, $\delta\in (0,1]$, $\alpha\in [\alpha_0,1)$, $\lambda \in \Sigma_\theta$, $f\in \L^2(\IR^d)$ and $F \in \L^2 (\IR^d ; \IC^{d })$ the solution $u\in \S^{\alpha,2}(\IR^d)$ to
    \begin{align}
    \label{eq: resolvent equation}
        \lambda u - \div^\alpha(\mu \nabla^\alpha u) = f + \div^\alpha(F) \quad \text{in }\IR^d
    \end{align}
    satisfies 
    \begin{align*}
        &|\lambda|\|u\|_{\L^2(B)}^2+\|\nabla^\alpha u\|_{\L^2(B)}^2\leq \frac{C}{|\lambda|}\|f\|_{\L^2(2B)}^2 \\
        &\quad + C  \sum\limits_{k=0}^\infty 2^{-\nu k} \bigg(\|F\|^2_{\L^2(B_{k+1})}+\delta \|\nabla^\alpha u\|^2_{\L^2(B_{k+1})}+\frac{\|u\|^2_{\L^2(B_{k+1})}}{\delta r^{2\alpha}}\bigg).
    \end{align*}
    The constant $C>0$ depends only on $d,\theta,\alpha_0,\nu$ and $\mu_\bullet,\mu^\bullet$.
\end{lemma}

\begin{proof}
    Fix $\varphi\in \C_c^\infty(2B)$ satisfying $0\leq \varphi\leq 1$, $\varphi=1$ on $B$ and $r\|\nabla \varphi \|_{\L^\infty}\lesssim 1$. By Lemma~\ref{lem: Leibniz}, we have $u\varphi^2\in \S^{\alpha,2}(\IR^d)$ and $\nabla^\alpha (u\varphi^2) = \varphi^2\nabla^\alpha u + R$, where $R=u \nabla^\alpha (\varphi^2) + \nabla^\alpha_{\NL}(u,\varphi^2)$. Testing the resolvent equation~\eqref{eq: resolvent equation} with $u\varphi^2$ leads to
    \begin{align}
    \label{eq: resolv 1}
        \lambda\int\limits_{\IR^d} |u\varphi|^2\,\d x + \int\limits_{\IR^d}\mu \nabla^\alpha u \cdot \overline{(\nabla^\alpha u)\varphi^2}\,\d x &= \int\limits_{\IR^d} f\overline{u\varphi^2} \,\d x-\int\limits_{\IR^d} F\cdot \overline{(\nabla^\alpha u) \varphi^2}\,\d x \notag\\
        &\quad - \int\limits_{\IR^d} (\mu \nabla^\alpha u+F)\cdot \overline{R}\,\d x.
    \end{align}
    Due to Assumption~\ref{Ass: Coefficients}, we notice that
    \begin{align*}
        \int\limits_{\IR^d}\mu \nabla^\alpha u \cdot \overline{(\nabla^\alpha u)\varphi^2}\,\d x   \in \overline{\Sigma_{\omega_0}}.
    \end{align*}
    Since $\lambda\in \Sigma_\theta$ and $\theta+\omega_0 <\pi$, there exists $C>0$ depending only on $\theta$, $\omega_0$ and $\mu_\bullet$, such that
    \begin{align*}
        |\lambda|\|u\varphi\|_{\L^2}^2+ \|\varphi \nabla^\alpha u\|_{\L^2}^2\leq C\bigg|\lambda\int\limits_{\IR^d} |u\varphi|^2\,\d x + \int\limits_{\IR^d}\mu \nabla^\alpha u \cdot \overline{(\nabla^\alpha u)\varphi^2}\,\d x\bigg|.
    \end{align*}
    Combining this estimate with the identity \eqref{eq: resolv 1} yields the bound
    \begin{align*}
        |\lambda|\|u\varphi\|_{\L^2}^2+ \|\varphi \nabla^\alpha u\|_{\L^2}^2&\leq C\int\limits_{\IR^d} |f| |u\varphi^2| \,\d x+C\int\limits_{\IR^d} |F| |\varphi^2\nabla^\alpha u|\,\d x\\
        &\quad + C\int\limits_{\IR^d} (|\mu \nabla^\alpha u|+|F|) |R| \,\d x.
    \end{align*}
    Applying Young's inequality to the first two terms on the right-hand side implies
    \begin{align*}
        C\int\limits_{\IR^d} |f| |u\varphi^2| \,\d x+C\int\limits_{\IR^d} |F| |\varphi^2\nabla^\alpha u|\,\d x&\leq \frac{2C^2}{|\lambda|}\int\limits_{\IR^d} |f\varphi|^2 \,\d x+ \frac{|\lambda|}{2}\int\limits_{\IR^d}|u\varphi|^2\,\d x\\
        &\quad +2C^2\int\limits_{\IR^d} |F\varphi|^2 \,\d x  +\frac{1}{2}\int\limits_{\IR^d}  |(\nabla^\alpha u) \varphi|^2\,\d x .
    \end{align*}
    Plugging this into the inequality above and rearranging the terms yields
    \begin{align}
    \begin{split}
    \label{eq: ident 2}
        |\lambda|\|u\varphi\|_{\L^2}^2+ \|\varphi \nabla^\alpha u\|_{\L^2}^2&\lesssim  \frac{1}{|\lambda|}\|f\varphi\|_{\L^2}^2 + \|F\varphi\|_{\L^2}^2+\int\limits_{\IR^d} (|\mu \nabla^\alpha u|+|F|) |R| \,\d x \\
        &\lesssim \frac{1}{|\lambda|}\|f\|_{\L^2(2B)}^2 + \|F\|_{\L^2(2B)}^2+\sum\limits_{k=0}^\infty  a_k\|R\|_{\L^2(C_k(B))},\raisetag{10pt}
    \end{split}
    \end{align}
    where we set $a_k \coloneqq \|\mu \nabla^\alpha u\|_{\L^2(C_k(B))} + \|F\|_{\L^2(C_k(B))}$ for $k\geq 0$ and used the properties of $\varphi$ and Hölder's inequality in the last step.
    Finally, we estimate the last term on the right-hand side. By Lemma~\ref{lem: auxiliary lem}, we have
    \begin{align*}
        &\|R\|_{\L^2(2B)}\lesssim  \frac{1}{r^\alpha}\sum\limits_{m=0}^\infty  2^{-(\frac{d}{2}+\alpha)m}\|u\|_{\L^2(C_m(B))},\\
        &\|R\|_{\L^2(C_k(B))}\lesssim\frac{2^{-(d+\alpha)k}}{r^\alpha}\|u\|_{\L^2(C_k(B))}+ \frac{2^{-(\frac{d}{2}+\alpha)k}}{r^\alpha}\|u\|_{\L^2(2B)} 
    \end{align*}
    for $k\geq 1$. Thus, we have
    \begin{align*}
        \sum\limits_{k=0}^\infty  a_k\|R\|_{\L^2(C_k(B))}
        & \lesssim\frac{a_0}{r^\alpha}\sum\limits_{m=0}^\infty  2^{-(\frac{d}{2}+\alpha)m}\|u\|_{\L^2(C_m(B))}\\
        &\quad +\sum\limits_{k=1}^\infty\bigg(\frac{2^{-(d+\alpha)k}a_k}{r^\alpha}\|u\|_{\L^2(C_k(B))}+ \frac{2^{-(\frac{d}{2}+\alpha)k}a_k}{r^\alpha}\|u\|_{\L^2(2B)}\bigg).
    \end{align*}
    For every fixed $\delta\in(0,1]$, we apply Young's inequality on the terms on the right-hand side to obtain
    \begin{align*}
        &\sum\limits_{k=0}^\infty  a_k\|R\|_{\L^2(C_k(B))}\\
        &\lesssim  \sum\limits_{m=0}^\infty \bigg( \delta 2^{(\nu-d-2\alpha)m} a_0^2  +\frac{2^{-\nu m}}{\delta r^{2\alpha}}\|u\|^2_{\L^2(C_m(B))} \bigg)\\
        &+\sum\limits_{k=1}^\infty\bigg(\delta2^{-(d+\alpha)k}a_k^2 + \frac{2^{-(d+\alpha)k}}{\delta r^{2\alpha}}\|u\|^2_{\L^2(C_k(B))} + \delta2^{-\nu k}a_k^2 + \frac{2^{(\nu-d-2\alpha)k}}{\delta r^{2\alpha}}\|u\|^2_{\L^2(2B)}\bigg).
    \end{align*}
    Since $\nu <d+\alpha_0$, we have $\sum_{k=0}^\infty 2^{(\nu-d-2\alpha)k}\leq \sum_{k=0}^\infty 2^{(\nu-d-2\alpha_0)k}<\infty.$
    Hence, we estimate
    \begin{align*}
        &\sum\limits_{k=0}^\infty  a_k\|R\|_{\L^2(C_k(B))}\\
        &\lesssim \delta a_0^2 + \sum\limits_{m=0}^\infty  \frac{2^{-\nu m}}{\delta r^{2\alpha}}\|u\|^2_{\L^2(C_m(B))}+\sum\limits_{k=1}^\infty\bigg(\delta2^{-\nu k}a_k^2 + \frac{2^{-(d+\alpha)k}}{\delta r^{2\alpha}}\|u\|^2_{\L^2(C_k(B))} \bigg) + \frac{\|u\|^2_{\L^2(2B)}}{\delta r^{2\alpha}}\\
        &\lesssim \delta \sum\limits_{k=0}^\infty2^{-\nu k}a_k^2 + \sum\limits_{k=0}^\infty \frac{2^{-\nu k}}{\delta r^{2\alpha}}\|u\|^2_{\L^2(C_k(B))}.
    \end{align*}
    Inserting this estimate into \eqref{eq: ident 2} and using $\varphi=1$ on $B$ and $C_k(B)\subset B_{k+1}$ proves the claim.
\end{proof}

\begin{theorem}[Nonlocal Caccioppoli inequality]
\label{thm: Caccioppoli}
    Let $\alpha_0\in(0,1)$. For all $\theta \in [0,\pi-\omega_0)$ and all $\nu \in (0, d+\alpha_0)$ there exists a constant $C >0$ such that for all balls $B= B(x_0,r)$, $\alpha\in [\alpha_0,1)$, $\lambda \in \Sigma_\theta$, $f\in \L^2(\IR^d)$ and $F \in \L^2 (\IR^d ; \IC^{d })$ the solution $u\in \S^{\alpha,2}(\IR^d)$ to \eqref{eq: resolvent equation}
    satisfies 
    \begin{align*}
        &|\lambda|\sum\limits_{\ell=0}^\infty 2^{-\nu\ell}\|u\|^2_{\L^2(B_\ell)}+\sum\limits_{\ell=0}^\infty 2^{-\nu\ell} \|\nabla^\alpha u\|^2_{\L^2(B_\ell)}\\
        &\leq   C\sum\limits_{\ell=0}^\infty 2^{-\nu \ell} \bigg(\frac{1}{|\lambda|}\|f\|_{\L^2(B_{\ell+1})}^2 +\|F\|_{\L^2(B_{\ell+1})}^2 +  \frac{1}{r^{2\alpha }}  \|u\|^2_{\L^2(B_{\ell+1})}\bigg).
    \end{align*}
    The constant $C>0$ depends only on $d,\alpha_0,\theta,\nu$ and $\mu_\bullet,\mu^\bullet$.
\end{theorem}

\begin{proof}
    Let $\nu' \coloneqq \frac{d + \alpha_0 + \nu}{2} \in (\nu, d + \alpha_0)$. We apply Lemma~\ref{lem: iteration setup} with $\nu'$ to the balls $B_\ell$, $\ell \in \IN_0$, whose radii are $2^\ell r \geq r$ and whose dilation satisfy $(B_\ell)_{k+1} = B_{k+\ell+1}$ to obtain 
    \begin{align*}
        &|\lambda|\|u\|_{\L^2(B_\ell)}^2+\|\nabla^\alpha u\|_{\L^2(B_\ell)}^2\leq \frac{C}{|\lambda|}\|f\|_{\L^2(B_{\ell+1})}^2\\
        &\quad +C  \sum\limits_{k=0}^\infty 2^{-\nu' k} \bigg(\|F\|^2_{\L^2(B_{k+\ell+1})}+\delta \|\nabla^\alpha u\|^2_{\L^2(B_{k+\ell+1})}+\frac{\|u\|^2_{\L^2(B_{k+\ell+1})}}{\delta r^{2\alpha}}\bigg).
    \end{align*}
    After the index shift $n = k + \ell$ we obtain
    \begin{align*}
        &|\lambda| \|u\|_{\L^2(B_\ell)}^2 + \|\nabla^\alpha u\|_{\L^2(B_\ell)}^2
        \leq \frac{C}{|\lambda|} \|f\|_{\L^2(B_{\ell+1})}^2 \\
        &\quad + C\sum_{n=\ell}^\infty 2^{-\nu'(n-\ell)} \Big( \|F\|_{\L^2(B_{n+1})}^2 + \delta \|\nabla^\alpha u\|_{\L^2(B_{n+1})}^2 + \frac{ \|u\|^2_{\L^2(B_{n+1})} }{\delta r^{2\alpha}}\Big).
    \end{align*}
    We multiply both sides by $2^{-\nu \ell}$ and sum over all $\ell \in \IN_0$. Since, for non-negative $(b_n)_n$,
    \begin{align*}
        \sum_{\ell=0}^\infty 2^{-\nu\ell}\sum_{n=\ell}^\infty 2^{-\nu'(n-\ell)} b_n = \sum_{n=0}^\infty 2^{-\nu'n} b_n \sum_{\ell=0}^n 2^{(\nu'-\nu)\ell} \leq C_1 \sum_{n=0}^\infty 2^{-\nu n} b_n
    \end{align*}
    with $C_1 = C_1(d,\alpha_0,\nu)>0$, we arrive at 
    \begin{align*}
        &|\lambda| \sum\limits_{\ell=0}^\infty 2^{-\nu \ell}\|u\|_{\L^2(B_\ell)}^2 + \sum\limits_{\ell=0}^\infty 2^{-\nu \ell}\|\nabla^\alpha u\|_{\L^2(B_\ell)}^2
        \leq \frac{C}{|\lambda|}\sum\limits_{\ell=0}^\infty 2^{-\nu \ell} \|f\|_{\L^2(B_{\ell+1})}^2 \\
        &\quad + CC_1\sum_{n=0}^\infty 2^{-\nu n} \Big( \|F\|_{\L^2(B_{n+1})}^2 + \delta \|\nabla^\alpha u\|_{\L^2(B_{n+1})}^2 + \frac{ \|u\|^2_{\L^2(B_{n+1})} }{\delta r^{2\alpha}}\Big).
    \end{align*}
    Using the identity
    \begin{align*}
        C C_1 \delta \sum_{n=0}^\infty 2^{-\nu n} \|\nabla^\alpha u\|_{\L^2(B_{n+1})}^2 =2^\nu C C_1 \delta \sum_{\ell=1}^\infty 2^{-\nu \ell} \|\nabla^\alpha u\|_{\L^2(B_\ell)}^2
    \end{align*}
    on the right-hand side, we can absorb the term with the fractional gradients into the left-hand side by choosing $\delta = \min\{(2^{\nu+1} C C_1)^{-1},2^{-1}\}$ to conclude
    \begin{align*}
        &|\lambda| \sum\limits_{\ell=0}^\infty 2^{-\nu \ell}\|u\|_{\L^2(B_\ell)}^2 + \sum\limits_{\ell=0}^\infty 2^{-\nu \ell}\|\nabla^\alpha u\|_{\L^2(B_\ell)}^2
        \\
        &\leq  C\sum\limits_{\ell=0}^\infty 2^{-\nu \ell} \bigg(\frac{1}{|\lambda|}\|f\|_{\L^2(B_{\ell+1})}^2 +\|F\|_{\L^2(B_{\ell+1})}^2 +  \frac{1}{r^{2\alpha }}  \|u\|^2_{\L^2(B_{\ell+1})}\bigg)
    \end{align*}
    for some constant $C>0$ depending only on $d,\alpha_0,\theta,\nu$ and $\mu_\bullet,\mu^\bullet$.
\end{proof}

\begin{theorem}[$\L^2$ off-diagonal estimates]
\label{thm: L2 off-diagonal estimates}
    Let $\alpha_0\in(0,1)$. For all $\theta \in [0,\pi-\omega_0)$ and all $\nu \in (0, d+\alpha_0)$ there exists a constant $C >0$ such that for all balls $B= B(x_0,r)$, $\alpha\in [\alpha_0,1)$, $\lambda \in \Sigma_\theta$, $f\in \L^2(\IR^d)$ and $F \in \L^2 (\IR^d ; \IC^{d })$ the solution $u\in \S^{\alpha,2}(\IR^d)$ to \eqref{eq: resolvent equation} satisfies 
    \begin{align*}
        &\sum\limits_{n=0}^\infty 2^{-\nu n}\| u\|_{\L^2(B_n)}^2+\frac{1}{|\lambda|}\sum\limits_{n=0}^\infty 2^{-\nu n}\|\nabla^\alpha u\|_{\L^2(B_n)}^2\\
        &\leq\frac{C}{|\lambda|}\sum\limits_{n=0}^{\infty}\sum\limits_{k=0}^n\bigg(\prod_{s=0}^{k-1}\frac{C}{|\lambda|(2^sr)^{2\alpha}}\bigg) 2^{-\nu (n-k)} \bigg(\frac{1}{|\lambda|}\|f\|_{\L^2(B_{n+1})}^2 +\|F\|_{\L^2(B_{n+1})}^2\bigg).
    \end{align*}
    The constant $C>0$ depends only on $d,\alpha_0,\theta,\nu$ and $\mu_\bullet,\mu^\bullet$.
\end{theorem}

\begin{proof}
    For $r>0$ and $B= B(x_0,r)$ set
    \begin{align*}
        S(r) = \sum\limits_{\ell=0}^\infty 2^{-\nu \ell}\|u\|_{\L^2(B_\ell)}^2 \quad \text{and}\quad T(r)=\sum\limits_{\ell=0}^\infty 2^{-\nu \ell} \bigg(\frac{1}{|\lambda|}\|f\|_{\L^2(B_{\ell+1})}^2 +\|F\|_{\L^2(B_{\ell+1})}^2\bigg).
    \end{align*}
    By Theorem~\ref{thm: Caccioppoli}, we have 
     \begin{align*}
        &S(r) +\frac{1}{|\lambda|}\sum\limits_{n=0}^\infty 2^{-\nu n}\|\nabla^\alpha u\|_{\L^2(B_n)}^2
        \leq   \frac{C}{|\lambda|}T(r)+\frac{C}{|\lambda|r^{2\alpha}}S(2r) .
    \end{align*}
    Applying this inequality successively with growing radii $2r,4r, \dots, 2^mr$ for $m\in\IN$ to the $S$-terms, we obtain
    \begin{align*}
        S(r) +\frac{1}{|\lambda|}\sum\limits_{n=0}^\infty 2^{-\nu n}\|\nabla^\alpha u\|_{\L^2(B_n)}^2
        &\leq   \frac{C}{|\lambda|}\sum\limits_{k=0}^{m-1}\bigg(\prod_{s=0}^{k-1}\frac{C}{|\lambda|(2^sr)^{2\alpha}}\bigg)T(2^kr)\\
        &\quad +\bigg(\prod_{s=0}^{m-1}\frac{C}{|\lambda|(2^sr)^{2\alpha}}\bigg)S(2^mr) .
    \end{align*}
    Since the prefactor $\prod_{s=0}^{m-1}\frac{C}{|\lambda|(2^sr)^{2\alpha}}\to 0$ as $m\to \infty$ and $S(2^mr) \lesssim \|u\|_{\L^2}^2<\infty$ for all $m\in\IN$, we obtain in the limit
    \begin{align*}
        S(r) +\frac{1}{|\lambda|}\sum\limits_{n=0}^\infty 2^{-\nu n}\|\nabla^\alpha u\|_{\L^2(B_n)}^2
        \leq   \frac{C}{|\lambda|}\sum\limits_{k=0}^{\infty}\bigg(\prod_{s=0}^{k-1}\frac{C}{|\lambda|(2^sr)^{2\alpha}}\bigg)T(2^kr).
    \end{align*}
    Finally, by the substitution $\ell+k=n$, we observe
    \begin{align*}
       T(2^kr) = \sum\limits_{n=k}^\infty 2^{-\nu (n-k)} \bigg(\frac{1}{|\lambda|}\|f\|_{\L^2(B_{n+1})}^2 +\|F\|_{\L^2(B_{n+1})}^2\bigg).
    \end{align*}
    Inserting this into the inequality above yields
    \begin{align*}
        &S(r) +\frac{1}{|\lambda|}\sum\limits_{n=0}^\infty 2^{-\nu n}\|\nabla^\alpha u\|_{\L^2(B_n)}^2\\
        &\leq   \frac{C}{|\lambda|}\sum\limits_{k=0}^{\infty}\sum\limits_{n=k}^\infty\bigg(\prod_{s=0}^{k-1}\frac{C}{|\lambda|(2^sr)^{2\alpha}}\bigg) 2^{-\nu (n-k)} \bigg(\frac{1}{|\lambda|}\|f\|_{\L^2(B_{n+1})}^2 +\|F\|_{\L^2(B_{n+1})}^2\bigg)\\
        &=\frac{C}{|\lambda|}\sum\limits_{n=0}^{\infty}\sum\limits_{k=0}^n\bigg(\prod_{s=0}^{k-1}\frac{C}{|\lambda|(2^sr)^{2\alpha}}\bigg) 2^{-\nu (n-k)} \bigg(\frac{1}{|\lambda|}\|f\|_{\L^2(B_{n+1})}^2 +\|F\|_{\L^2(B_{n+1})}^2\bigg),
    \end{align*}
    which proves the claim.
\end{proof}

\begin{corollary}
\label{cor: simplyfied coefficient estimates}
    Let $\alpha_0\in(0,1)$. For all $\nu \in (0,d+\alpha_0)$ there exists a constant $C>0$ such that for all $t>0$, $\alpha\in [\alpha_0,1)$, $f\in \L^2(\IR^d)$, $F\in\L^2(\IR^d;\IC^d)$ and balls $B=B(x_0,r)$ with $r\geq t$ we have
    \begin{align*}
        &\|t^\alpha \nabla^\alpha(1+t^{2\alpha}L)^{-1}f\|^2_{\L^2(B)}+\|t^\alpha(1+t^{2\alpha}\cL)^{-1}\div^\alpha F\|^2_{\L^2(B)}\\
        &\leq C\sum\limits_{n=0}^\infty 2^{-\nu n} \Big(\|f\|_{\L^2(C_n(B))}^2+ \|F\|_{\L^2(C_n(B))}^2\Big).
    \end{align*}
    The constant $C>0$ depends only on $d,\alpha_0,\nu$ and $\mu_\bullet,\mu^\bullet$.
\end{corollary}

The proof follows the same line of reasoning as in \cite[Cor.~2.6]{Haardt_Tolksdorf} and will therefore be omitted.

\begin{remark}
\label{rem: cubes for balls}
    One can replace balls by cubes in the statements of Theorem~\ref{thm: L2 off-diagonal estimates} and Corollary~\ref{cor: simplyfied coefficient estimates}, since the nonlocal Caccioppoli inequality in Theorem~\ref{thm: Caccioppoli} also holds for cubes by applying Remark~\ref{rem: cubes instead of balls}.
\end{remark}

\section{Proof of the main result}
\label{sec: proof}
This section is devoted to proving our main result, Theorem~\ref{thm: Kato property}. In particular, we develop a second-order approach, which follows closely the algorithm of \cite{Haardt_Tolksdorf}; see also \cites{AHLMcT} and \cites{ISem27}.
Throughout this section, we fix $\alpha_0\in(0,1)$, $\alpha\in[\alpha_0,1)$ and $\mu\in \L^\infty(\IR^d;\IC^{d\times d})$ satisfying Assumption~\ref{Ass: Coefficients}. Moreover, unless stated otherwise, all implicit constants suppressed by $\lesssim$ depend only on $d$, $\alpha_0$, $\mu_\bullet$ and $\mu^\bullet$. In particular, they are independent of $\alpha\in[\alpha_0,1)$.

\subsection{Reduction to a square function estimate}
\label{Sec: Reduction to a square function estimate}

It is well known that the Kato property has a deep connection to certain square function estimates. In particular, the Kato property in Theorem~\ref{thm: Kato property} with constants depending only on $d,\alpha_0$ and $\mu_\bullet,\mu^\bullet$ is equivalent to the square function estimates 
\begin{align}
\label{eq:SFE start solution Kato}
        \int\limits_0^\infty \|t^\alpha L(1+t^{2\alpha}L)^{-1}u\|^2_{\L^2}\,\frac{\d t}{t} &\leq C \|\nabla^\alpha u\|_{\L^2}^2,\\
        \int\limits_0^\infty \|t^\alpha L^*(1+t^{2\alpha}L^*)^{-1}u\|^2_{\L^2}\,\frac{\d t}{t} &\leq C \|\nabla^\alpha u\|_{\L^2}^2\nonumber 
\end{align}
for every $u\in\S^{\alpha,2}(\IR^d)$, where $C>0$ depends only on $d,\alpha_0$ and $\mu_\bullet,\mu^\bullet$; see \cite[Prop.~12.7]{ISem27} with the substitution $s=t^{\alpha}$. Recall that $L^*$ is again a fractional operator associated to the coefficients $\mu^*$ that satisfy Assumption~\ref{Ass: Coefficients} with the same ellipticity constants, see Section~\ref{sec: second-order operator}.  Hence, for Theorem~\ref{thm: Kato property}, it suffices to show the first square function estimate in~\eqref{eq:SFE start solution Kato}.

As a first step, we rewrite the square function estimate~\eqref{eq:SFE start solution Kato} in terms of $\nabla^\alpha u$. For this purpose, we use the identity from Lemma~\ref{lem: identity for cL} to obtain
\begin{align}
\label{eq: rewrite}
    t^\alpha L(1+t^{2\alpha}L)^{-1}u = t^\alpha (1+t^{2\alpha}\cL)^{-1}\cL u = -t^\alpha (1+t^{2\alpha}\cL)^{-1} \div^\alpha(\mu \nabla^\alpha u).
\end{align}
Employing the uniform $\L^2$-bounds from Proposition~\ref{prop: L2 resolvent bounds}, we can define the right-hand side as a bounded operator on vector-valued functions.
\begin{definition}
\label{def: Theta_t}
For $t>0$ we define the operator
\begin{align*}
     \Theta_t :\L^2(\IR^d;\IC^{d}) \to \L^2(\IR^d), \quad  \Theta_t F \coloneqq  -t^\alpha (1+t^{2\alpha} \cL)^{-1}\div^\alpha  (\mu F).
\end{align*}
\end{definition}

Hence, using \eqref{eq: rewrite}, we see that the square function estimate~\eqref{eq:SFE start solution Kato} is equivalent to
\begin{align}
\label{eq: SFE for Theta}
    \int_0^\infty \|\Theta_t(\nabla^\alpha u)\|_{\L^2}^2 \, \frac{\d t}{t} \lesssim \|\nabla^\alpha u\|_{\L^2}^2 \qquad (u \in \S^{\alpha,2}(\IR^d)).
\end{align}

\subsection{Principal part approximation}

The idea of the principal part approximation is to approximate the operator $\Theta_t$ in an averaged sense in order to reduce the square function estimate \eqref{eq: SFE for Theta} to a Carleson measure estimate. For this purpose, we next introduce the following dyadic structure: let $\square = \bigcup_{j\in\IZ} \square_{2^j}$ where $\square_{2^j} = \{2^j(k+(0,1]^d): k\in\IZ^d\}$. A dyadic cube $Q\in \square_{2^j}$ has side-length $\ell(Q) = 2^j$ and measure $|Q| = 2^{jd}$. Moreover define $\square_{t} = \square_{2^j}$ for $2^{j-1}<t\leq 2^j$. Then the dyadic averaging operator $\cA_t$ on $\L^2(\IR^d)$ is defined via
\begin{align*}
    \cA_t u(x) \coloneqq \fint_Q u(y)\, \d y 
\end{align*}
for every $x\in \IR^d$ and $t>0$, where $Q$ is the unique dyadic cube belonging to $\square_t$ and containing the point $x$. Note that $\cA_t^2 = \cA_t$. Furthermore, when a scalar operator, such as $\cA_t$, is applied to a vector-valued function, we agree from now on to apply the operator to each component separately.

As a first step of the principal part approximation, we want to extend the domain of $\Theta_t$ from Definition~\ref{def: Theta_t} to bounded vector fields. For this purpose, recall that $(t^{\alpha}(1+t^{2\alpha}\cL)^{-1}\div^\alpha)_{t>0}$ satisfies the nonlocal off-diagonal bounds from Corollary~\ref{cor: simplyfied coefficient estimates}, which also hold for balls replaced by cubes, see Remark~\ref{rem: cubes for balls}. Since $\Theta_t$ arises from that operator family applied to the bounded multiplication operator $\mu$, these estimates are transferred to $\Theta_t$ and will be recorded in the next lemma. 
\begin{lemma}
\label{lem: off-diag for Theta_t}
    For all $\nu \in (0,d+\alpha_0)$ there exists a constant $C>0$ depending only on $d,\alpha_0, \nu$ and $\mu_\bullet,\mu^\bullet$ such that for all $t>0$, $F\in\L^2(\IR^d;\IC^d)$ and cubes $Q$ with side-length $\ell(Q)\geq t$ we have
    \begin{align*}
        &\|\Theta_tF\|_{\L^2(Q)}^2\leq C\sum\limits_{n=0}^\infty 2^{-\nu n} \|F\|_{\L^2(C_n(Q))}^2.
    \end{align*}
\end{lemma}

\begin{lemma}
\label{lem: extension of Theta to Linf}
    Let $b\in \L^\infty(\IR^d;\IC^{d})$. Then for every $t>0$ and cube $Q\subset\IR^d$ with sides parallel to the coordinate axes the limit
    \begin{align*}
        \Theta_t b \coloneqq \lim_{j\to \infty}\Theta_t(\ind_{2^jQ}b)
    \end{align*}
    exists in $\L_\loc^2(\IR^d)$ and is independent of the underlying cube $Q$.
\end{lemma}

The proof is similar to \cite[Prop.~4.2]{Haardt_Tolksdorf} and only needs the off-diagonal estimates from Lemma~\ref{lem: off-diag for Theta_t}. The remaining details are omitted for brevity. 

\begin{definition}
    Identify the standard unit vectors $e_1,\dots, e_{d} \in \IC^{d}$ with the respective constant functions on $\IR^d$. For $t>0$ we define the \textit{principal part} of $\Theta_t$ by
    \begin{align*}
        [\gamma_t]_{j} \coloneqq \Theta_t(e_j) \in \Lloc^2(\IR^d)
    \end{align*}
    for $j=1,\dots,d$.
\end{definition}

To state the main reduction step of this subsection, we briefly recapitulate the notion of a Carleson measure.

\begin{definition}
\label{def: Carleson measure}
    A Borel measure $\Lambda$ on $\IR^{d+1}_+$ is called a \textit{Carleson measure} if there exists a constant $C>0$ such that
    \begin{align*}
        \Lambda (R(Q))\leq C|Q|\qquad (Q\in \square),
    \end{align*}
    where $R(Q) = (0,\ell(Q)]\times Q$. 
\end{definition}

\begin{proposition}
\label{prop:reduction to Carleson}
    If
    \begin{align}
    \label{eq: principal part approximation}
        \int\limits_0^\infty \|(\Theta_t - \gamma_t\cdot \cA_t)(\nabla^\alpha u)\|_{\L^2}^2 \, \frac{\d t}{t} \lesssim \|\nabla^\alpha u\|_{\L^2}^2 
    \end{align}
    holds for every $u \in \S^{\alpha,2}(\IR^d)$ and $|\gamma_t(x)|^2 \, \frac{\d x \, \d t}{t}$ is a Carleson measure, then~\eqref{eq: SFE for Theta} is valid for all $u \in \S^{\alpha,2}(\IR^d)$.
\end{proposition}

Its proof is an application of Carleson's lemma, following the same reasoning as \cite[Prop.~4.5]{Haardt_Tolksdorf}. The rest of this subsection is devoted to verifying \eqref{eq: principal part approximation}. This will be done in the same fashion as in \cite[Sec.~4]{Haardt_Tolksdorf}. For this purpose we make use of the following identity: for a fixed cube $Q\in \square_t$ and $F\in \L^2(\IR^d;\IC^{d})$ we have for $x\in Q$
\begin{align*}
    (\gamma_t\cdot \cA_t)F(x) = \sum_{j=1}^{d} \Theta_t(e_j)(x) (F_j)_Q  = \Theta_t\Big(\sum_{j=1}^{d}(F_j)_Q e_j \Big) (x)=  \Theta_t((F)_Q)(x).
\end{align*}
Furthermore, by Lemma~\ref{lem: extension of Theta to Linf}, the following holds on $Q$
\begin{align}
\label{eq: rep of prinicple part}
    \begin{split}
        (\gamma_t\cdot \cA_t )F &= \lim\limits_{k\to \infty} \Theta_t(\mathbf{1}_{2^kQ} (F)_Q)\\
    &= \lim\limits_{k\to \infty}\sum\limits_{n=0}^{k-1} \Theta_t(\mathbf{1}_{C_n(Q)} (F)_Q) = \sum\limits_{n=0}^\infty \Theta_t(\mathbf{1}_{C_n(Q)} (F)_Q).
    \end{split}
\end{align}

\begin{lemma}
\label{lem: uniform bounds of principal part}
    For all $t>0$ and $F\in \L^2(\IR^d;\IC^{d})$ we have
\begin{align*}
    \|(\gamma_t\cdot \cA_t)F\|_{\L^2} \lesssim \|F\|_{\L^2}.
\end{align*}
\end{lemma}

\begin{proof}
    For fixed $t>0$ it is enough to show the bound
    \begin{align*}
        \|(\gamma_t\cdot \cA_t)F\|_{\L^2(Q)}^2\lesssim \|F\|^2_{\L^2(Q)}
    \end{align*}
    for all $Q\in \square_t$. Since $\ell(Q) = 2^j\ge t$, using the identity~\eqref{eq: rep of prinicple part} and Lemma~\ref{lem: off-diag for Theta_t} yields
    \begin{align*}
         \|(\gamma_t\cdot \cA_t)F\|_{\L^2(Q)} &\leq  \sum\limits_{n=0}^\infty\| \Theta_t(\mathbf{1}_{C_n(Q)} (F)_Q)\|_{\L^2(Q)}\lesssim  \sum\limits_{n=0}^\infty 2^{-\frac{\nu}{2} n} \|(F)_Q\|_{\L^2(C_n(Q))}
    \end{align*}
    for all $0<\nu<d+\alpha_0$. Due to Hölder's inequality, we have
    \begin{align*}
        \|(F)_Q\|_{\L^2(C_n(Q))} \leq \frac{|C_n(Q)|^\frac{1}{2}}{|Q|}\int\limits_{Q}|F|\,\d x\leq\frac{|C_n(Q)|^\frac{1}{2}}{|Q|^\frac{1}{2}} \|F\|_{\L^2(Q)} \simeq 2^{\frac{d}{2}n}\|F\|_{\L^2(Q)}.
    \end{align*}
    Hence, we derive 
    \begin{align*}
        \|(\gamma_t\cdot \cA_t)F\|_{\L^2(Q)} \lesssim \sum\limits_{n=0}^\infty 2^{(\frac{d}{2}-\frac{\nu}{2} )n}\|F\|_{\L^2(Q)} \simeq \|F\|_{\L^2(Q)},
    \end{align*}
    where we chose $\nu= d+\frac{\alpha_0}{2}>d$ in the last step.
\end{proof}

The following lemma presents an approximation of $\Theta_t$ in an averaged sense in space but pointwise in time $t$. Note that such a control by a fractional Sobolev norm was also used in the resolution of the parabolic Kato problem and for the generalized Stokes operator, see \cite[Sec.~7.3]{Auscher_Egert_Nystrom} and \cite[Lem.~4.7]{Haardt_Tolksdorf}.

\begin{lemma}
\label{lem:uniform bound principal part approximation}
For all $t>0$ we have
\begin{align*}
    \|(\Theta_t - \gamma_t \cdot  \cA_t)F\|_{\L^2} \lesssim t^\beta \| F\|_{\dot\W^{\beta,2}} \qquad (F \in \S^{\beta,2}(\IR^d;\IC^{d}))
\end{align*}
for every $\beta\in (0,\frac{\alpha_0}{2})$. The implicit constant depends only on $d,\alpha_0,\beta$ and $\mu_\bullet, \mu^\bullet$.
\end{lemma}

\begin{proof}
Fix $t>0$. We argue as in \cite[Lem.~4.7]{Haardt_Tolksdorf} in two steps. First we work on a fixed dyadic cube, and then sum over a partition of $\IR^d$.\\

\noindent \textbf{Step~1: Estimate on a fixed cube $Q \in \square_t$.}
Fix a cube $Q\in\square_t$. Since $\Theta_t$ is a bounded operator on $\L^2$, we have
\begin{align*}
    \Theta_t F =  \sum\limits_{\ell=0}^{\infty}\Theta_t \bigl(\ind_{C_\ell(Q)}F\bigr).
\end{align*}
By subtracting the identity \eqref{eq: rep of prinicple part} and taking the $\L^2(Q)$-norm, we obtain
\begin{align*}
    \|(\Theta_t - \gamma_t\cdot  \cA_t)F\|_{\L^2(Q)} 
    & \leq  \sum\limits_{\ell=0}^{\infty} \bigl\|\Theta_t \bigl(\mathbf{1}_{C_\ell(Q)}(F - (F)_Q) \bigr)\bigr\|_{\L^2(Q)}.
\end{align*}
Because $\ell(Q)\geq t$, we can apply Lemma~\ref{lem: off-diag for Theta_t} to derive
\begin{align}
\label{eq: theta 1}
     \begin{split}
        \|(\Theta_t - \gamma_t \cdot \cA_t)F\|_{\L^2(Q)} 
        &\lesssim    \sum\limits_{\ell=0}^{\infty}\bigg(\sum\limits_{n=0}^{\infty} 2^{-\nu n} \bigl\|\mathbf{1}_{C_\ell(Q)}(F - (F)_Q) \bigr\|^2_{\L^2(C_n(Q))}\bigg)^\frac{1}{2}\\
        &\leq   \sum\limits_{\ell=0}^{\infty}2^{-\frac{\nu}{2} \ell} \|F - (F)_Q \|_{\L^2(2^{\ell+1}Q)}.
     \end{split}
\end{align}
Next, we introduce a telescoping sum to estimate
\begin{align*}
    \|F - (F)_Q \|_{\L^2(2^{\ell+1}Q)}& \leq \|F - (F)_{2^{\ell+1}Q} \|_{\L^2(2^{\ell+1}Q)} + \sum\limits_{j=1}^{\ell+1}\|(F)_{2^{j}Q}  - (F)_{2^{j-1}Q} \|_{\L^2(2^{\ell+1}Q)}\\
    & = \|F - (F)_{2^{\ell+1}Q} \|_{\L^2(2^{\ell+1}Q)} + \sum\limits_{j=1}^{\ell+1}|2^{\ell+1}Q|^\frac{1}{2}|(F)_{2^{j}Q}  - (F)_{2^{j-1}Q}|\\
    &=  \|F - (F)_{2^{\ell+1}Q} \|_{\L^2(2^{\ell+1}Q)} + \sum\limits_{j=1}^{\ell+1}\frac{|2^{\ell+1}Q|^\frac{1}{2}}{|2^{j-1}Q|}\bigg| \int\limits_{2^{j-1}Q} F- (F)_{2^jQ}\,\d x\bigg|\\
    & \lesssim \sum\limits_{j=1}^{\ell+1}2^{\frac{d}{2}(\ell+1 -j)} \|F- (F)_{2^jQ}\|_{\L^2(2^{j}Q)},
\end{align*}
where we used Hölder's inequality in the last step. Now, we apply the fractional Poincar\'e inequality from \cite[Lem.~2.2]{fracSobolev} componentwise for $\beta\in (0,\frac{\alpha_0}{2})$ to get
\begin{align}
\label{eq: frac Sobolev-Poincare}
    \|F - (F)_{2^{j}Q} \|_{\L^2(2^{j}Q)}
    \lesssim |2^jQ|^\frac{\beta}{d} \|F\|_{\dot\W^{\beta,2}(2^jQ)} 
    \lesssim t^\beta 2^{\beta j} \|F\|_{\dot\W^{\beta,2}(2^jQ)}
\end{align}
for all $j\in \IN$ and implicit constants depending only on $d$ and $\beta$. This leads to the bound
\begin{align*}
    \|F - (F)_Q \|_{\L^2(2^{\ell+1}Q)}& \lesssim    t^\beta \sum\limits_{j=1}^{\ell+1} 2^{\frac{d}{2}(\ell+1-j) +\beta j} \|F\|_{\dot\W^{\beta,2} (2^jQ)}.
\end{align*}
Inserting this estimate into \eqref{eq: theta 1}, squaring both sides and applying Hölder's inequality, we obtain for $0<\nu'<\nu$ that
\begin{align}
\label{eq: single Q}
    \begin{split}
        \|(\Theta_t - \gamma_t\cdot  \cA_t)F\|^2_{\L^2(Q)}& \lesssim t^{2\beta }\sum\limits_{\ell=0}^{\infty}2^{-\nu' \ell}\bigg( \sum\limits_{j=1}^{\ell+1} 2^{\frac{d}{2}(\ell+1-j) +\beta j} \|F\|_{\dot\W^{\beta,2} (2^jQ)}\bigg)^2 \\
        &\lesssim t^{2\beta } \sum\limits_{\ell=0}^{\infty}(\ell+1)2^{-\nu' \ell}\sum\limits_{j=1}^{\ell+1} 2^{d(\ell+1-j) +2\beta j}  \|F\|^2_{\dot\W^{\beta,2} (2^jQ)}.
    \end{split}
\end{align}

\noindent \textbf{Step~2: Conclusion}.
Summing over all $Q \in \square_t$ in \eqref{eq: single Q} gives us
\begin{align*}
    \|(\Theta_t - \gamma_t \cdot \cA_t)F\|_{\L^2}^2\lesssim t^{2\beta } \sum\limits_{\ell=0}^{\infty}(\ell+1)2^{-\nu' \ell}\sum\limits_{j=1}^{\ell+1} 2^{d(\ell+1-j) +2\beta j} \sum_{Q \in \square_t} \|F\|^2_{\dot\W^{\beta,2} (2^jQ)}.
\end{align*}
By a counting argument, we observe that
\begin{align*}
     \sum_{Q \in \square_t}\|F\|^2_{\dot\W^{\beta,2} (2^jQ)} 
     \leq    \int\limits_{\IR^d}\int\limits_{\IR^d} \biggl(\sum_{Q \in \square_t} \mathbf{1}_{2^jQ}(x)\biggr)\frac{|F(x)-F(y)|^2}{|x-y|^{d+2\beta }}\,\d x\d y \lesssim 2^{dj}  \|F\|^2_{\dot\W^{\beta,2} (\IR^d)}
\end{align*}
for all $j\in\IN$. Thus, we conclude
\begin{align*}
    \|(\Theta_t - \gamma_t\cdot \cA_t)F\|_{\L^2}^2
    &\lesssim t^{2\beta } \sum\limits_{\ell=0}^{\infty}(\ell+1)2^{-\nu' \ell}\sum\limits_{j=1}^{\ell+1} 2^{d(\ell+1) +2\beta j}\|F\|^2_{\dot\W^{\beta,2} (\IR^d)} \\
    &\lesssim t^{2\beta }\|F\|^2_{\dot\W^{\beta,2} (\IR^d)}   \sum\limits_{\ell=0}^{\infty}(\ell+1)2^{(d+2\beta-\nu' )\ell} \\
    &\lesssim t^{2\beta }\|F\|^2_{\dot\W^{\beta,2} (\IR^d)},
\end{align*}
where the series converges for every $\beta \in (0,\frac{\alpha_0}{2})$ since we can choose $d+2\beta<\nu'<d+\alpha_0$.
\end{proof}

To prove~\eqref{eq: principal part approximation}, we will first smooth out the ``harsh'' approximation to get a control via quadratic estimates of some well-behaved operator and afterwards remove the smoothing again to gain back the full estimate. For this purpose, we introduce the following smoothing operators.

\begin{definition}
\label{def: P_t and Q_t}
    For $t>0$ we define the following operators on $\L^2(\IR^d)$:
    \begin{align*}
        P_t = (1+t^{2\alpha}(-\Delta)^\alpha)^{-1}\quad \text{and}\quad Q_t = t^{\alpha}\nabla^\alpha (1+t^{2\alpha}(-\Delta)^\alpha)^{-1}.
    \end{align*}
\end{definition}
As a reminder, the reader should keep in mind that we agreed to use the same notation for extensions of operators on vector fields by acting componentwise.

\begin{lemma}
\label{lem: properties of P_t and Q_t}
Let $(P_t)_{t>0}$ and $(Q_t)_{t>0}$ be as in Definition~\ref{def: P_t and Q_t}. Then we have the following assertions.
    \begin{enumerate}
        \item\label{item: 1} The operator families $(P_t)_{t>0}$ and $(Q_t)_{t>0}$ are uniformly $\L^2$-bounded. Moreover, for every $u\in\S^{\alpha,2}(\IR^d)$, we have $t^\alpha P_t\nabla^\alpha u = Q_t u$.
        
        \item \label{item: 2} We have the adjoint square function estimate
        \begin{align*}
            \int\limits_0^\infty \|Q^*_t \nabla^\alpha u\|^2_{\L^2}\,\frac{\d t}{t}=\frac{1}{2\alpha} \|\nabla^\alpha u\|^2_{\L^2} \qquad (u\in \S^{\alpha,2}(\IR^d)).
        \end{align*}

        \item\label{item: 3} The following Calderón reproducing formula holds
        \begin{align*}
            2\alpha\int\limits_0^\infty  Q_tQ_t^* \nabla^\alpha u\,\frac{\d t}{t}= \nabla^\alpha u \qquad (u\in \S^{\alpha,2}(\IR^d)).
        \end{align*}
    \end{enumerate}
\end{lemma}

\begin{proof}
We prove the three assertions separately.\\

\noindent \textbf{Step 1: Show \eqref{item: 1}}.
By Proposition~\ref{prop: L2 resolvent bounds} with $\mu = \Id_{\IC^{d\times d}}$ being the identity matrix, we observe that
\begin{align*}
    \|P_tu\|_{\L^2} +\|Q_tu\|_{\L^2} \lesssim \|u\|_{\L^2}
\end{align*}
for some implicit constant depending only on $d$. Furthermore, for all $u\in\S^{\alpha,2}(\IR^d)$ we have the identity $t^\alpha P_t\nabla^\alpha u = Q_t u$ by taking the Fourier transform, see Proposition~\ref{prop: fourier symbols}.\\

\noindent\textbf{Step 2: Show \eqref{item: 2}}.
Notice that $Q_t^* = -t^\alpha (1+t^{2\alpha}(-\Delta)^\alpha)^{-1}\div^\alpha$, where $(-\Delta)^\alpha=-\div^\alpha \nabla^\alpha$ has to be understood in the sense of Lemma~\ref{lem: identity for cL}. Hence, we compute
\begin{align}
\label{eq: Qt*}
    Q_t^*\nabla^\alpha u = -t^\alpha (1+t^{2\alpha}(-\Delta)^\alpha)^{-1}\div^\alpha(\nabla^\alpha u) = t^\alpha (-\Delta)^\alpha (1+t^{2\alpha}(-\Delta)^\alpha)^{-1}u
\end{align}
for every $u\in\S^{\alpha,2}(\IR^d)$. Define $\psi_\alpha(r)\coloneqq r^\alpha(1+r^{2\alpha})^{-1}$ for $r>0$. Then, we get with Plancherel's and Fubini's theorems
\begin{align*}
    \int\limits_0^\infty \|Q^*_t \nabla^\alpha u\|^2_{\L^2}\,\frac{\d t}{t} &= \int\limits_0^\infty \|t^\alpha (-\Delta)^\alpha (1+t^{2\alpha}(-\Delta)^\alpha)^{-1}u\|^2_{\L^2}\,\frac{\d t}{t}\\
    &= \int\limits_0^\infty \int\limits_{\IR^d} \big||\xi|^\alpha\psi_\alpha(t|\xi|)\cF(u)(\xi)\big|^2\,\frac{\d \xi\d t}{t}\\
    &=  \int\limits_{\IR^d} \bigg(\int\limits_0^\infty |\psi_\alpha(t|\xi|)|^2\frac{\d t}{t}\bigg) \big||\xi|^\alpha \cF(u)(\xi)\big|^2\,\d \xi.
\end{align*}
Using the substitution $t|\xi| \mapsto t$, Proposition~\ref{prop: fourier symbols} with $|\xi|^\alpha=\big|\xi|\xi|^{\alpha-1}\big|$ and Plancherel's theorem leaves us with
\begin{align*}
    \int\limits_0^\infty \|Q^*_t \nabla^\alpha u\|^2_{\L^2}\,\frac{\d t}{t}= \int\limits_{\IR^d} \bigg(\int\limits_0^\infty |\psi_\alpha(t)|^2\frac{\d t}{t}\bigg) \big|\xi|\xi|^{\alpha-1} \cF(u)(\xi)\big|^2\,\d \xi = \frac{1}{2\alpha}\|\nabla^\alpha u\|_{\L^2}^2,
\end{align*}
where we used the integral identity
\begin{align*}
    \int\limits_0^\infty |\psi_\alpha(t)|^2\frac{\d t}{t}=\int\limits_0^\infty t^{2\alpha}(1+t^{2\alpha})^{-2}\frac{\d t}{t} = \frac{1}{2\alpha}\int\limits_0^\infty s(1+s)^{-2}\frac{\d s}{s} = \frac{1}{2\alpha}.
\end{align*}

\noindent\textbf{Step 3: Show \eqref{item: 3}}.
Applying $Q_t$ to \eqref{eq: Qt*} yields 
\begin{align*}
    Q_tQ_t^*\nabla^\alpha u =t^{2\alpha} (-\Delta)^\alpha (1+t^{2\alpha}(-\Delta)^\alpha)^{-2}\nabla^\alpha u,
\end{align*}
for $u\in\S^{\alpha,2}(\IR^d)$ by commuting the corresponding symbols in Fourier space with Proposition~\ref{prop: fourier symbols}. Hence, \eqref{item: 3} follows, for example, from the Calderón reproducing formula of the negative Laplacian $(-\Delta)$, see \cite[Thm.~6.16]{ISem27}.
\end{proof}

To prove \eqref{eq: principal part approximation}, we decompose
\begin{align}
\label{eq: split of ppa}
    \Theta_t - \gamma_t\cdot \cA_t =  \Theta_t(1- P_t)+(\Theta_t - \gamma_t\cdot \cA_t)P_t -  \gamma_t\cdot \cA_t(1- P_t)
\end{align}
and treat the three terms separately.

\begin{proposition}
\label{prop: ppa 1}
    For every $u\in\S^{\alpha,2}(\IR^d)$ we have
    \begin{align*}
        \int\limits_0^\infty \|\Theta_t(1- P_t) \nabla^\alpha u\|_{\L^2}^2\,\frac{\d t}{t}\lesssim \|\nabla^\alpha u\|_{\L^2}^2.
    \end{align*}
\end{proposition}

\begin{proof}
By Definition~\ref{def: Theta_t}, \eqref{eq: rewrite} and \eqref{eq: Qt*}, we have for every $u\in\S^{\alpha,2}(\IR^d)$ the identity
\begin{align*}
    \Theta_t(1- P_t) \nabla^\alpha u = \Theta_t\nabla^\alpha (1- P_t)  u
    &=-t^\alpha (1+t^{2\alpha} \cL)^{-1}\div^\alpha (\mu \nabla^\alpha (1- P_t)  u)\\
    &=t^\alpha L(1+t^{2\alpha}L)^{-1} t^{2\alpha}(-\Delta)^\alpha (1+t^{2\alpha} (-\Delta)^\alpha)^{-1}u\\
    &=(1-(1+t^{2\alpha}L)^{-1}) Q_t^*\nabla^\alpha u.
\end{align*}
Combining this with Proposition~\ref{prop: L2 resolvent bounds} yields the estimate
\begin{align*}
    \|\Theta_t(1- P_t) \nabla^\alpha u\|_{\L^2}\lesssim \|Q_t^*\nabla^\alpha u\|_{\L^2},
\end{align*}
which concludes the claim by Lemma~\ref{lem: properties of P_t and Q_t}.
\end{proof}

\begin{proposition}[Smoothed principal part approximation]
\label{prop: ppa 2}
    For all $F\in \L^2(\IR^d;\IC^d)$ we have
    \begin{align*}
        \int\limits_0^\infty \|(\Theta_t - \gamma_t\cdot \cA_t)P_t F\|_{\L^2}^2\,\frac{\d t}{t}\lesssim \|F\|_{\L^2}^2.
    \end{align*}
\end{proposition}

\begin{proof}
    Set $\beta=\frac{\alpha_0}{4}$. By Lemma~\ref{lem:uniform bound principal part approximation} and Proposition~\ref{prop: homogeneous norms}, we have
    \begin{align*}
        \|(\Theta_t - \gamma_t\cdot \cA_t)P_t F\|_{\L^2} \lesssim t^{\beta}\| P_t F\|_{\dot\W^{\beta,2}}\simeq \| t^{\beta}(-\Delta)^\frac{\beta}{2} (1+t^{2\alpha}(-\Delta)^{\alpha})^{-1}F \|_{\L^2}.
    \end{align*}
    Finally, using Plancherel and Fubini's theorems yields
    \begin{align*}
        \int\limits_0^\infty \|(\Theta_t - \gamma_t\cdot \cA_t)P_t F\|_{\L^2}^2\,\frac{\d t}{t} &\lesssim \int\limits_0^\infty \| t^{\beta}(-\Delta)^\frac{\beta}{2} (1+t^{2\alpha}(-\Delta)^{\alpha})^{-1}F \|^2_{\L^2}\,\frac{\d t}{t}\\
        &=\int\limits_0^\infty\int\limits_{\IR^d} |t^{\beta}|\xi|^\beta (1+t^{2\alpha}|\xi|^{2\alpha})^{-1}\cF(F)(\xi)|^2\,\frac{\d \xi\d t}{t}\\
        &=\int\limits_{\IR^d} \bigg(\int\limits_0^\infty (t|\xi|)^{2\beta} (1+(t|\xi|)^{2\alpha})^{-2}\,\frac{\d t}{t}\bigg)|\cF(F)(\xi)|^2\,\d \xi\\
        &=\bigg(\int\limits_0^\infty t^{2\beta} (1+t^{2\alpha})^{-2}\,\frac{\d t}{t}\bigg)\|F\|^2_{\L^2}.
    \end{align*}
    Using the estimate of the integral
    \begin{align*}
        \int\limits_0^\infty t^{2\beta} (1+t^{2\alpha})^{-2}\,\frac{\d t}{t} \leq \int\limits_0^1 t^{2\beta} \,\frac{\d t}{t}+\int\limits_1^\infty t^{-2\alpha}\,\frac{\d t}{t} =  \frac{1}{2\beta} + \frac{1}{2\alpha }\leq \frac{5}{2\alpha_0},
    \end{align*}
    we obtain the claim with a constant that depends only on $d,\alpha_0$ and $\mu_\bullet, \mu^\bullet$.
\end{proof}

To handle the last term of the splitting~\eqref{eq: split of ppa}, we need the following auxiliary lemma, which is also crucial for the $T(b)$ argument. We refer to \cite[Lem.~13.6]{ISem27} for a local version.

\begin{lemma}
\label{lem:palmen}
Let $\alpha_0\in (0,1)$. There exists $C>0$ such that for all $0<\gamma<\alpha_0$, $Q \in \square$, $\alpha\in[\alpha_0,1)$ and all $u \in \S^{\alpha,2}(\IR^d)$ we have
\begin{align*}
    \bigg|\fint_Q \nabla^\alpha u \, \d x\bigg|^2 &\leq  \frac{C}{ \ell(Q)^{\frac{2\alpha}{3}}}\bigg(\sum\limits_{k=0}^\infty 2^{-\gamma k}\fint_{2^{k+1}Q}  |u|^2 \, \d x \bigg)^\frac{1}{3}\bigg(\fint_Q |\nabla^\alpha u|^2\, \d x\bigg)^\frac{2}{3}.
\end{align*}
The constant $C>0$ depends only on $d$ and $\alpha_0$.
\end{lemma}

\begin{proof}
Let $a \in (0,1)$ be chosen later. Fix $u \in \S^{\alpha,2}(\IR^d)$, $Q\in\square$ with side-length $\ell= \ell(Q)$, and a function $\eta \in \C_c^\infty(Q)$ with $0\leq \eta\leq 1$ and $\eta = 1$ on $(1-a)Q$ that satisfies $a \ell \|\nabla \eta\|_{\L^\infty} \leq C$ for some constant $C$ that only depends on $d$. Applying Definition~\ref{def: S space} to $F=\eta e_j$ for every $1\leq j\leq d$, we get $\int \eta \nabla^\alpha u \,\d x = -\int (\nabla^\alpha \eta) u\,\d x$. Thus, with the splitting $\nabla^\alpha u = \eta \nabla^\alpha u + (1-\eta) \nabla^\alpha  u$ we obtain
\begin{align*}
    \fint_Q \nabla^\alpha u \, \d x 
    = \fint_Q \eta \nabla^\alpha u \, \d x  + \fint_Q (1-\eta) \nabla^\alpha u \, \d x = - \frac{1}{|Q|}\int_{\IR^d} (\nabla^\alpha \eta)  u \, \d x  + \fint_Q (1-\eta) \nabla^\alpha u \, \d x.
\end{align*}
Set $E \coloneqq Q \setminus (1-a)Q$. With H\"older's inequality we derive
\begin{align*}
    \bigg|\fint_Q \nabla^\alpha u \, \d x\bigg|
    &\leq\frac{1}{|Q|}\int\limits_{\IR^d} |\nabla^\alpha  \eta | |u| \, \d x + C \fint_Q \mathbf{1}_E |\nabla^\alpha u|\, \d x \notag\\
    &\leq \frac{1}{|Q|}\int\limits_{\IR^d} |\nabla^\alpha  \eta| |u| \, \d x  + C \biggl(\frac{|E|}{|Q|}\biggr)^{\frac{1}{2}}\bigg(\fint_Q |\nabla^\alpha u|^2\, \d x\bigg)^{\frac{1}{2}}.
\end{align*}
For the first term on the right-hand side, we apply Lemma~\ref{lem: auxiliary lem}, see also Remark~\ref{rem: cubes instead of balls}, to obtain
\begin{align*}
    \int\limits_{\IR^d} |\nabla^\alpha  \eta | |u| \, \d x= \sum\limits_{k=0}^\infty\, \int\limits_{C_k(Q)} |\nabla^\alpha  \eta | |u| \, \d x  \leq  C\sum\limits_{k=0}^\infty \frac{2^{-(d+\alpha_0) k}}{(a\ell)^\alpha}\int\limits_{2^{k+1}Q}  |u| \, \d x ,
\end{align*}
which leads to the estimate
\begin{align*}
    \bigg|\fint_Q \nabla^\alpha u \, \d x\bigg| &\leq \frac{C}{(a\ell)^\alpha}\sum\limits_{k=0}^\infty 2^{-\alpha_0 k}\fint_{2^{k+1}Q}  |u| \, \d x + C \biggl(\frac{|E|}{|Q|}\biggr)^{\frac{1}{2}}\bigg(\fint_Q |\nabla^\alpha u|^2\, \d x\bigg)^{\frac{1}{2}}.
\end{align*}
Squaring both sides, estimating $\frac{|E|}{|Q|} = (1-(1-a)^d) \leq a d\leq a^\alpha d$ by Bernoulli's inequality, and using Hölder's inequality for $0<\gamma<\alpha_0$, we arrive at
\begin{align*}
    \bigg|\fint_Q \nabla^\alpha u \, \d x\bigg|^2 &\leq C\Big(\frac{X}{a^{2\alpha}\ell^{2\alpha }} + a^\alpha Y \Big)
\end{align*}
for some constant $C>0$ depending on $d$ and $\alpha_0$, where
\begin{align*}
    X\coloneqq \sum\limits_{k=0}^\infty 2^{-\gamma k}\fint_{2^{k+1}Q}  |u|^2 \, \d x  \quad \text{and} \quad Y\coloneqq\fint_Q |\nabla^\alpha u|^2\, \d x. 
\end{align*}
Optimizing the inequality in the parameter $a$ leads to the value
\begin{align*}
    a_0 \coloneqq  \Big( \frac{2 X}{\ell^{2\alpha} Y}\Big)^\frac{1}{3\alpha}.
\end{align*}
If $a_0<1$, then we can choose $a=a_0$ so that the inequality becomes
\begin{align*}
    \bigg|\fint_Q \nabla^\alpha u \, \d x\bigg|^2 &\leq \frac{C}{\ell^{\frac{2\alpha}{3}}} X^\frac{1}{3}Y^\frac{2}{3},
\end{align*}
which proves the claim. In the case $a_0\geq 1$, this condition reads
\begin{align*}
    \fint_Q |\nabla^\alpha  u|^2\, \d x\leq\frac{2}{\ell^{2\alpha}}\sum\limits_{k=0}^\infty 2^{-\gamma k}\fint_{2^{k+1}Q}  |u|^2 \, \d x .
\end{align*}
Hence, using Jensen's inequality, we directly calculate
\begin{align*}
    \bigg|\fint_Q \nabla^\alpha  u \, \d x\bigg|^2 &\leq \bigg(\fint_Q |\nabla^\alpha u|^2 \, \d x\bigg)^\frac{1}{3}\bigg(\fint_Q |\nabla^\alpha u|^2 \, \d x\bigg)^\frac{2}{3}\\
    &\leq \frac{2^\frac{1}{3}}{ \ell^{\frac{2\alpha}{3}}}\bigg(\sum\limits_{k=0}^\infty 2^{-\gamma k}\fint_{2^{k+1}Q}  |u|^2 \, \d x \bigg)^\frac{1}{3}\bigg(\fint_Q |\nabla^\alpha u|^2\, \d x\bigg)^\frac{2}{3},
\end{align*}
so that the claim follows once again.
\end{proof}

\begin{proposition}
\label{prop: ppa 3}
    For every $u\in\S^{\alpha,2}(\IR^d)$ we have
    \begin{align}
    \label{eq: assertion 1}
        \int\limits_0^\infty \|\cA_t(1  - P_t)\nabla^\alpha u \|_{\L^2}^2\,\frac{\d t}{t} \lesssim \|\nabla^\alpha u \|_{\L^2}^2.
    \end{align}
    Furthermore, we have the same bounds for $\gamma_t\cdot \cA_t(1-P_t)$.
\end{proposition}

\begin{proof}
    We follow the same line of reasoning as in \cite[Prop.~13.13]{ISem27}, see also \cite[Prop.~5.7]{Axelsson_Keith_McIntosh}.
    Since $\cA^2_t = \cA_t$ and $\gamma_t\cdot \cA_t$ is uniformly $\L^2$-bounded by Lemma~\ref{lem: uniform bounds of principal part}, it suffices to prove~\eqref{eq: assertion 1}.\\

    \noindent \textbf{Step 1: Reduction to a Schur-type bound.}
    In this step, we show that \eqref{eq: assertion 1} follows from
    \begin{align}
    \label{eq1:third term estimate}
        \|\cA_s (1- P_s) Q_tf \|_{\L^2} \lesssim \zeta(st^{-1}) \|f \|_{\L^2}\quad ( s,t>0, f\in\L^{2}(\IR^d)),
    \end{align}
    where $\zeta(\tau) \coloneqq \min (\tau^\alpha, \tau^{-\alpha}+ \tau^{-\frac{\alpha}{3}})\in \L^1((0,\infty); \frac{\d \tau}{\tau})$.
    In fact, by Lemma~\ref{lem: properties of P_t and Q_t} we have
    \begin{align*}
        \cA_s (1- P_s)\nabla^\alpha u = 2\alpha \int\limits_{0}^\infty \cA_s (1- P_s)Q_tQ^*_t \nabla^\alpha u\,\frac{\d t}{t}.
    \end{align*}
    Combining this representation with \eqref{eq1:third term estimate} yields the estimate
    \begin{align*}
        \int\limits_0^\infty \|\cA_s (1- P_s) \nabla^\alpha u \|_{\L^2}^2 \, \frac{\d s}{s} 
        &\lesssim \int\limits_0^\infty \biggl( \int\limits_0^\infty \|\cA_s (1- P_s) Q_t Q_t^* \nabla^\alpha u \|_{\L^2} \, \frac{\d t}{t} \biggr)^2 \, \frac{\d s}{s} \\
        &\lesssim \int\limits_0^\infty \biggl( \int\limits_0^\infty \zeta(st^{-1}) \|Q_t^*  \nabla^\alpha u \|_{\L^2} \, \frac{\d t}{t} \biggr)^2 \, \frac{\d s}{s}.
    \end{align*}
    Finally, using Hölder's inequality, the substitution $\tau = st^{-1}$ and Fubini's theorem gives us
    \begin{align*}
        \int\limits_0^\infty \|\cA_s (1- P_s)\nabla^\alpha u \|_{\L^2}^2 \, \frac{\d s}{s}
        &\lesssim \int\limits_0^\infty \biggl( \int\limits_0^\infty \zeta(st^{-1})\, \frac{\d t}{t} \biggr) \biggl( \int\limits_0^\infty \zeta(st^{-1}) \|Q_t^*\nabla^\alpha u \|^2_{\L^2} \, \frac{\d t}{t} \biggr)\, \frac{\d s}{s} \\
        &=\int\limits_0^\infty \biggl( \int\limits_0^\infty \zeta(\tau)\, \frac{\d \tau}{\tau} \biggr) \biggl( \int\limits_0^\infty \zeta(st^{-1}) \|Q_t^*\nabla^\alpha u \|^2_{\L^2} \, \frac{\d t}{t} \biggr)\, \frac{\d s}{s} \\
        &=\bigg(\int\limits_0^\infty \zeta(\tau) \, \frac{\d \tau}{\tau}\bigg)^2 \biggl( \int\limits_0^\infty \|Q_t^*  \nabla^\alpha u \|_{\L^2}^2 \, \frac{\d t}{t} \biggr).
    \end{align*}
    Using the quadratic estimates from Lemma~\ref{lem: properties of P_t and Q_t} and the bound
    \begin{align*}
        \int\limits_0^\infty \zeta(\tau) \, \frac{\d \tau}{\tau}\leq \int\limits_0^1 \tau^\alpha \, \frac{\d \tau}{\tau}+\int\limits_1^\infty \tau^{-\alpha}+\tau^{-\frac{\alpha}{3}} \, \frac{\d \tau}{\tau} \leq \frac{5}{\alpha}\leq \frac{5}{\alpha_0},
    \end{align*}
    we conclude \eqref{eq: assertion 1} from the inequality above.\\

    \noindent \textbf{Step 2: Show \eqref{eq1:third term estimate}.}
    Using the Fourier description from Proposition~\ref{prop: fourier symbols}, we obtain the identities
    \begin{align*}
        P_sQ_t = \cF^{-1}\Big(\frac{1}{1+s^{2\alpha}|\xi|^{2\alpha}}\frac{it^\alpha \xi}{|\xi|^{1-\alpha}(1+t^{2\alpha}|\xi|^{2\alpha})}\cF\Big) =  \frac{t^\alpha}{s^\alpha}P_tQ_s
    \end{align*}
    and
    \begin{align*}
         (1- P_s)Q_t =  \cF^{-1}\Big(\frac{s^{2\alpha}|\xi|^{2\alpha}}{1+s^{2\alpha}|\xi|^{2\alpha}}\frac{it^\alpha \xi}{|\xi|^{1-\alpha}(1+t^{2\alpha}|\xi|^{2\alpha})}\cF\Big)=\frac{s^{\alpha}}{t^{\alpha}} (1 - P_t)Q_s.
    \end{align*}
    When $s\leq t$ we obtain \eqref{eq1:third term estimate} by uniform $\L^2$-bounds as follows
    \begin{align*}
        \|\cA_s (1- P_s) Q_tf\|_{\L^2} \lesssim\|(1- P_s) Q_tf\|_{\L^2} =\frac{s^{\alpha}}{t^{\alpha}}\|(1- P_t) Q_sf\|_{\L^2} \lesssim \frac{s^{\alpha}}{t^{\alpha}}\|f\|_{\L^2}.
    \end{align*}
    When $t\leq s$, we first estimate
    \begin{align}
    \label{eq2:third term estimate}
        \|\cA_s (1- P_s) Q_tf\|_{\L^2} &\lesssim \|\cA_s Q_t f\|_{\L^2} + \|P_s Q_t f\|_{\L^2} \notag\\
        &= \|\cA_s Q_t f\|_{\L^2} + \frac{t^\alpha}{s^\alpha} \|P_t Q_s f\|_{\L^2} \\
        &\lesssim \|\cA_s Q_t f\|_{\L^2} + \frac{t^\alpha}{s^\alpha} \|f\|_{\L^2}.\notag
    \end{align}
    To control the term $\|\cA_s Q_t f\|_{\L^2} $ on the right-hand side, we first decompose the $\L^2$-norm by dyadic cubes to get
    \begin{align*}
        \|\cA_s Q_t f\|_{\L^2}^2 = \sum_{Q \in \square_s} |Q| \biggl| \fint_Q Q_t f \, \d x \biggr|^2  = \sum_{Q \in \square_s} |Q| \biggl| \fint_Q t^\alpha \nabla^\alpha  P_t f \, \d x \biggr|^2.
    \end{align*}
    For $\gamma=\frac{\alpha_0}{2}$ we apply Lemma~\ref{lem:palmen} and  Hölder's inequality to get
    \begin{align*}
        \|\cA_s Q_t f\|_{\L^2}^2&\lesssim \sum_{Q \in \square_s} \frac{|Q|}{s^{\frac{2\alpha}{3}}} \bigg(\sum\limits_{k=0}^\infty 2^{-\frac{\alpha_0}{2} k}\fint_{2^{k+1}Q}  |t^\alpha P_t f|^2 \, \d x\bigg)^\frac{1}{3}\bigg(\fint_Q |t^\alpha \nabla^\alpha P_t f|^2\, \d x\bigg)^\frac{2}{3} \\
        &\leq \frac{1}{s^{\frac{2\alpha}{3}}}  \bigg(\sum_{Q \in \square_s}\sum\limits_{k=0}^\infty 2^{-(d+\frac{\alpha_0}{2}) k}\int\limits_{2^{k+1}Q}  |t^\alpha P_t f|^2 \, \d x\bigg)^\frac{1}{3}\bigg(\sum_{Q \in \square_s}\int\limits_Q |t^\alpha \nabla^\alpha P_t f|^2\, \d x\bigg)^\frac{2}{3}.
    \end{align*}
    With the observation $\sum_{Q\in\square_s}\ind_{2^kQ}(x)\lesssim 2^{kd}$ for all $k\in\IN_0$, we conclude
    \begin{align*}
         \|\cA_s Q_t f\|_{\L^2}^2&\lesssim \frac{1}{s^{\frac{2\alpha}{3}}}  \bigg(\sum\limits_{k=0}^\infty 2^{-\frac{\alpha_0}{2} k}\int\limits_{\IR^d}  |t^\alpha P_t f|^2 \, \d x\bigg)^\frac{1}{3}\bigg(\int\limits_{\IR^d} |t^\alpha \nabla^\alpha P_t f|^2\, \d x\bigg)^\frac{2}{3}\\
         &\lesssim \Big(\frac{t}{s}\Big)^\frac{2\alpha}{3}\|P_tf\|_{\L^2}^\frac{2}{3}\|Q_tf\|_{\L^2}^\frac{4}{3}\\
         &\lesssim \Big(\frac{t}{s}\Big)^\frac{2\alpha}{3}\|f\|_{\L^2}^2.
    \end{align*}
    This together with \eqref{eq2:third term estimate} completes the proof of \eqref{eq1:third term estimate} in the case $t\leq s$.
\end{proof}

Combining Propositions~\ref{prop: ppa 1}, \ref{prop: ppa 2} and \ref{prop: ppa 3} with the splitting \eqref{eq: split of ppa} establishes the full principal part approximation.

\begin{proposition}[Principal part approximation]
\label{prop: Full principal part approximation}
    For all $u\in \S^{\alpha,2}(\IR^d)$ we have
    \begin{align*}
        \int\limits_0^\infty \|(\Theta_t- \gamma_t\cdot \cA_t)\nabla^\alpha u\|_{\L^2}^2\,\frac{\d t}{t}\lesssim \|\nabla^\alpha u\|_{\L^2}^2.
    \end{align*}
\end{proposition}

\subsection{The Carleson measure argument}
\label{sec: Tb argument}
By Proposition~\ref{prop: Full principal part approximation} and Proposition~\ref{prop:reduction to Carleson}, it remains to show that $|\gamma_t(x)|^2\frac{\d x\d t}{t}$ is a Carleson measure. Following the usual procedure, we will first establish a directional Carleson measure estimate, where we restrict the values of the principal part to a cone in $\IC^d$.

\begin{definition}
\label{defi:cone in Cn}
Let $\eps>0$. For $\xi \in \IC^{d}$ with $|\xi| = 1$, define the open cone with central axis $\xi$ by
\begin{align*}
    \Gamma_\xi^\eps \coloneqq \bigl\{ w\in \IC^{d}\setminus \{0\} :\bigl|\tfrac{w}{|w|} - \xi\bigr| < \eps \bigr\}.
\end{align*}
\end{definition}

\begin{proposition}[Directional Carleson measure estimate]
\label{prop:key estimate Kato}
There is a choice of $\eps>0$, depending only on $d$, $\alpha_0$, $\mu_{\bullet}$ and $\mu^{\bullet}$, such that for every $\xi \in \IC^{d}$ with $\lvert \xi \rvert = 1$ we have
\begin{align*}
    \iint_{R(Q)} |\gamma_t(x)|^2 \, \ind_{\Gamma_\xi^\eps}(\gamma_t(x)) \, \frac{\d x \, \d t}{t} \lesssim |Q| \qquad (Q \in \square).
\end{align*}
\end{proposition}

Taking the directional Carleson measure estimate from Proposition~\ref{prop:key estimate Kato} for granted, we are ready to prove Theorem~\ref{thm: Kato property}.

\begin{proof}[Proof of Theorem~\ref{thm: Kato property}]
Let $\eps>0$ be as in Proposition~\ref{prop:key estimate Kato}. Then we cover $\IC^{d} \setminus \{ 0 \}$ by finitely many cones $\Gamma_{\xi_1}^\eps, \ldots, \Gamma_{\xi_N}^\eps$, where $N$ depends only on $d$ and $\eps$. The estimate from Proposition~\ref{prop:key estimate Kato} then yields
\begin{align*}
    \iint_{R(Q)} |\gamma_t(x)|^2 \, \frac{\d x \, \d t}{t}
    \leq \iint_{R(Q)} \sum_{j=1}^N |\gamma_t(x)|^2 \, \ind_{\Gamma_{\xi_j}^\eps}(\gamma_t(x)) \, \frac{\d x \, \d t}{t}
    \lesssim N |Q| 
\end{align*}
for every cube $Q \in \square$ as required. Thus, $|\gamma_t(x)|^2 \frac{\d x \d t}{t}$ defines a Carleson measure by Definition~\ref{def: Carleson measure}. Combining this with Propositions~\ref{prop: Full principal part approximation} and~\ref{prop:reduction to Carleson}, we obtain~\eqref{eq: SFE for Theta}, which proves the claim; see Section~\ref{Sec: Reduction to a square function estimate}.
\end{proof}

The remainder of this subsection is devoted to the proof of Proposition~\ref{prop:key estimate Kato}, which is based on the construction of a $T(b)$-test function. More precisely, for given $Q \in \square$ and $\xi \in \IC^{d}$ with $| \xi | = 1$, we strive for a construction of functions $b_Q^{\eps}$ such that
\begin{align*}
    |w|\lesssim \big|w\cdot (\cA_t b_Q^{\eps})(x) \big| \qquad (w \in \Gamma_{\xi}^{\eps})
\end{align*}
for a hopefully large set of $(t,x) \in R(Q)$. The test function $b_Q^{\eps}$ should have the following properties:
\begin{enumerate}
    \item It should belong to the range of $\nabla^\alpha$ to make the principal part approximation accessible.

    \item It should satisfy the $\L^2$-bound $\|b_Q^\eps\|_{\L^2}\lesssim |Q|^\frac{1}{2}$.
    
    \item It should point roughly in the direction of $\overline{\xi} \in \IC^{d}$ on $Q$.
\end{enumerate}
This will be achieved by the ansatz
\begin{align*}
    b_Q^{\eps} \sim \nabla^\alpha (1 + \eps^{2\alpha} \ell (Q)^{2\alpha} L)^{-1} \Phi.
\end{align*}
By construction, this function satisfies the first point. Moreover, since $(1 + \eps^{2\alpha} \ell(Q)^{2\alpha} L)^{-1}$ is an approximation of the identity for small $\eps$, properties (2) and (3) hold roughly if $\Phi$ is chosen such that $\|\nabla^\alpha\Phi\|_{\L^2}\lesssim|Q|^\frac{1}{2}$ and $\nabla^\alpha \Phi = \overline{\xi}$ on $Q$. As explained in the introduction, the elliptic construction of $\Phi$ as an affine function times a cut-off function is not admissible in our fractional situation. Thus, we instead use the factorization $\nabla^\alpha = \nabla I_{1 - \alpha}$. More precisely, composing with the fractional Laplacian restores the local gradient
\begin{align*}
    \nabla^\alpha (-\Delta)^{\frac{1-\alpha}{2}}
    = \nabla I_{1-\alpha} (-\Delta)^{\frac{1-\alpha}{2}} = \nabla ,
\end{align*}
see Proposition~\ref{prop: local correspondence}.
Therefore, we choose $\Phi = (-\Delta)^{\frac{1-\alpha}{2}} h$, where $h$ is the elliptic ansatz, so that $\nabla^\alpha \Phi = \nabla h$ has the required properties. This will be made rigorous in the proof of the following proposition.

\begin{proposition}
\label{prop:the Tb test functions}
There exists $\eps_0 \in (0,1]$ depending only on $d,\alpha_0, \mu_\bullet,\mu^\bullet$ such that for all $\xi\in\IC^{d}$ with $|\xi|=1$, all $ \eps \in(0, \eps_0]$ and all cubes $Q \in \square$ there exists $b_Q^\eps \in \L^2 (\IR^d ; \IC^{d})$ with the following properties:
\begin{enumerate}[label=(B\arabic*)]
    \item \label{item1: Tb test fct} $\|b_Q^\eps\|_{\L^2} \lesssim |Q|^{\frac{1}{2}}$,
    \item \label{item2: Tb test fct} $\displaystyle \Re \biggl(\xi \cdot \fint_Q b_Q^\eps \, \d x\biggr) \geq 1$,
    \item \label{item3: Tb test fct} $\displaystyle \iint_{R(Q)} |(\gamma_t(x)\cdot \cA_t) b_Q^\eps(x)|^2 \, \frac{\d x \, \d t}{t} \lesssim \frac{|Q|}{\eps^{2}}$.
\end{enumerate}
\end{proposition}

\begin{proof}
We fix $\xi\in \IC^{d}$ with $|\xi|=1$, $Q \in \square$ and abbreviate by $x_Q\in\IR^d$ and $\ell \coloneqq \ell(Q)$ its center and side-length. Similarly, we simplify notation by omitting $\eps$ and $Q$ when constructing the test function $b = b_Q^\eps$. Fix $\eta \in \C_c^\infty(\IR^d)$ with $\supp \eta \subset 2Q$, $\eta \equiv 1$ on $Q$ and $\|\eta\|_{\L^\infty}+\ell \|\nabla \eta\|_{\L^\infty}\lesssim 1$. Define 
\begin{align*}
    h (x)\coloneqq \bar \xi \cdot (x-x_Q)\eta(x)
    \in \C_c^\infty(\IR^d),
\end{align*}
and set $\Phi \coloneqq (-\Delta)^{\frac{1-\alpha}{2}}h $, which belongs to $ \S^{\alpha,2}(\IR^d)$ due to Proposition~\ref{prop: local correspondence}. Finally, we define our test function
\begin{align}
\label{eq2:the Tb test function}
    b \coloneqq  2\nabla^\alpha  (1+\eps^{2\alpha} \ell^{2\alpha}L )^{-1} \Phi.
\end{align}
By construction, we have with \eqref{eq: rewrite}
\begin{align}
\label{eq3:the Tb test function}
\begin{split}
    \frac{1}{2}b- \nabla^\alpha  \Phi
    &= \nabla^\alpha  \bigl((1+ \eps^{2\alpha} \ell^{2\alpha} L)^{-1} - 1 \bigr) \Phi \\
    &= -\nabla^\alpha  \bigl(\eps^{2\alpha} \ell^{2\alpha} L(1+ \eps^{2\alpha} \ell^{2\alpha} L)^{-1} \bigr) \Phi \\
    &=\eps^{2\alpha} \ell^{2\alpha} \nabla^\alpha(1+ \eps^{2\alpha} \ell^{2\alpha} \cL)^{-1} \div^\alpha(\mu \nabla^\alpha\Phi).
\end{split}
\end{align}
We are now ready to prove the properties~\ref{item1: Tb test fct}, \ref{item2: Tb test fct} and \ref{item3: Tb test fct}.\\

\noindent \textbf{Step 1: Show \ref{item1: Tb test fct}}.
Using Proposition~\ref{prop: local correspondence}, we calculate
\begin{align}
\label{eq: locla Phi}
    (\nabla^\alpha\Phi)(x) =\nabla^\alpha (-\Delta)^{\frac{1-\alpha}{2}}h(x)
    = \nabla h(x)
    =\bar \xi \eta(x)+  \bar \xi\cdot (x-x_Q) \nabla \eta(x).
\end{align}
Since $|x-x_Q|\leq \sqrt{d}\ell$ for every $x\in 2Q$ and $\ell\|\nabla\eta\|_{\L^\infty}\lesssim 1$, we obtain the pointwise bound
\begin{align*}
    |(\nabla^\alpha \Phi)(x)|\leq  |\bar \xi \eta(x)|+  |\bar \xi\cdot (x-x_Q) \nabla \eta(x)|\lesssim\Big( 1+\frac{\sqrt{d}\ell}{\ell} \Big)\ind_{2Q}(x),
\end{align*}
and subsequently
\begin{align}
\label{eq4:the Tb test function}
    \| \nabla^\alpha \Phi\|_{\L^2} \lesssim\|\ind_{2Q}\|_{\L^2}\lesssim|Q|^\frac{1}{2}.
\end{align}
Combining this estimate with the identity~\eqref{eq3:the Tb test function} and the uniform $\L^2$-bounds from Proposition~\ref{prop: L2 resolvent bounds}, we get
\begin{align*}
    \| b- 2\nabla^\alpha  \Phi\|_{\L^2}
    \lesssim \| \mu \nabla^\alpha \Phi\|_{\L^2}
    \lesssim |Q|^\frac{1}{2},
\end{align*}
which concludes \ref{item1: Tb test fct}.\\

\noindent \textbf{Step 2: Show \ref{item2: Tb test fct}}.
We see from \eqref{eq: locla Phi} that $\nabla^\alpha \Phi = \bar \xi$ on $Q$. Hence, by~\eqref{eq3:the Tb test function}, we have
\begin{align}
\label{eq5:the Tb test function}
    \fint_Q b - 2\overline{\xi} \, \d x= \fint_Q b - 2\nabla^\alpha  \Phi \, \d x=  2\eps^{2\alpha} \ell^{2\alpha} \fint_Q \nabla^\alpha u \, \d x,
\end{align}
where we set $u \coloneqq (1+\eps^{2\alpha} \ell^{2\alpha} \cL)^{-1}\div^\alpha(\mu\nabla^\alpha \Phi)$. The uniform $\L^2$-bounds from Proposition~\ref{prop: L2 resolvent bounds} in combination with \eqref{eq4:the Tb test function} give us
\begin{align*}
    \|u\|_{\L^2} \lesssim \frac{|Q|^{\frac{1}{2}} }{\eps^\alpha\ell^\alpha}\quad \text{and} \quad
    \|\nabla^\alpha  u\|_{\L^2} \lesssim \frac{|Q|^{\frac{1}{2}}}{\eps^{2\alpha} \ell^{2\alpha}},
\end{align*}
so that estimating the average on the right-hand side of \eqref{eq5:the Tb test function} by means of Lemma~\ref{lem:palmen} with $\gamma=\frac{\alpha_0}{2}$ leads us to
\begin{align*}
    \biggl|\fint_Q b - 2\overline{\xi} \, \d x \biggr|^2
    &\leq \frac{4C\eps^{4\alpha} \ell^{4\alpha}}{\ell^\frac{2\alpha}{3}}\bigg(\sum\limits_{k=0}^\infty 2^{-\frac{\alpha_0}{2} k}\fint_{2^{k+1}Q}  |u|^2 \, \d x\bigg)^\frac{1}{3} \bigg(\fint_Q |\nabla^\alpha u|^2\, \d x\bigg)^\frac{2}{3}\\
    &\leq \frac{4C\eps^{4\alpha} \ell^{4\alpha}}{ \ell^{\frac{2\alpha}{3}}|Q|}\|u\|_{\L^2}^{\frac{2}{3}}\|\nabla^\alpha u\|_{\L^2}^\frac{4}{3}\\
    &\leq C\eps^{\frac{2\alpha}{3}},
\end{align*}
where $C$ varies in every step. Since $|\xi|=1$, we conclude that
\begin{align*}
    \Re \biggl(\xi \cdot \fint_Q b \, \d x\biggr) 
    = 2 + \Re \biggl(\xi \cdot \fint_Q b - 2\overline{\xi} \, \d x \biggr) 
    &\geq  2 -\bigg| \fint_Q b - 2\overline{\xi} \, \d x \bigg| \\
    &\geq 2 - \sqrt{C \eps^{\frac{2\alpha}{3}} }\\
    &\geq 2 - \sqrt{C \eps^{\frac{2\alpha_0}{3}}}.
\end{align*}
Now, \ref{item2: Tb test fct} follows upon taking $\eps$ small enough.\\

\noindent \textbf{Step 3: Show \ref{item3: Tb test fct}}.
Recall from \eqref{eq2:the Tb test function} that $b\in \cR(\nabla^\alpha )$. Hence, we can apply the principal part approximation from Proposition~\ref{prop: Full principal part approximation} and use~\ref{item1: Tb test fct} to bound
\begin{align}
\label{eq6:the Tb test function}
\frac{1}{2} \iint_{R(Q)} |(\gamma_t(x)\cdot \cA_t) b(x)|^2 \, \frac{\d x  \d t}{t}
&\leq \int\limits_0^{\ell} \|(\gamma_t \cdot \cA_t - \Theta_t)b\|_{\L^2}^2  + \|\Theta_t b\|_{\L^2}^2 \, \frac{\d t}{t}\notag \\
&\lesssim\|b\|_{\L^2}^2  + \int\limits_0^{\ell} \|\Theta_t b\|_{\L^2}^2 \, \frac{\d t}{t} \\
&\lesssim|Q| +  \int\limits_0^{\ell} \|\Theta_t b\|_{\L^2}^2 \, \frac{\d t}{t}.\notag
\end{align}
By definition of $\Theta_t$ and $b$, we have
\begin{align*}
\Theta_t b 
&= -2t^\alpha(1+t^{2\alpha}\cL)^{-1}\div^\alpha(\mu \nabla^\alpha  (1+\eps^{2\alpha} \ell^{2\alpha}L )^{-1} \Phi)\\
&= 2t^\alpha(1+t^{2\alpha}\cL)^{-1}\cL (1+\eps^{2\alpha} \ell^{2\alpha}L )^{-1} \Phi\\
&= -t^\alpha(1+t^{2\alpha}L)^{-1} (1+\eps^{2\alpha} \ell^{2\alpha}\cL )^{-1}\div^\alpha(2 \mu\nabla^\alpha \Phi).
\end{align*}
Employing the uniform $\L^2$-bounds of Proposition~\ref{prop: L2 resolvent bounds} once again and \eqref{eq4:the Tb test function}, we obtain
\begin{align*}
    \|\Theta_t b\|_{\L^2}^2
    \lesssim \frac{t^{2\alpha}}{ \eps^{2\alpha} \ell^{2\alpha}} \|2\mu\nabla^\alpha \Phi\|_{\L^2}^2
    \lesssim \frac{ t^{2\alpha}}{ \eps^{2\alpha} \ell^{2\alpha}} |Q|.
\end{align*}
Finally, integration in $t$ leaves us with 
\begin{align*}
    \int\limits_0^{\ell} \|\Theta_t b\|_{\L^2}^2 \, \frac{\d t}{t} \lesssim \int\limits_0^{\ell} \frac{  t^{2\alpha}}{ \eps^{2\alpha} \ell^{2\alpha}} |Q| \, \frac{\d t}{t}
    \leq \frac{|Q|}{ 2\alpha \eps^{2\alpha}} \leq \frac{|Q|}{2\alpha_0\varepsilon^{2}},
\end{align*}
which we insert in \eqref{eq6:the Tb test function} to conclude \ref{item3: Tb test fct} provided $\eps \leq 1$.
\end{proof}

\begin{proof}[Proof of Proposition~\ref{prop:key estimate Kato}]
Fix $\varepsilon=\varepsilon_0$ from Proposition~\ref{prop:the Tb test functions}. Given the test functions  $b_Q^\eps$, the directional Carleson measure estimate follows from a stopping time argument with a John--Nirenberg lemma for Carleson measures, as in \cite{AHLMcT}, in the form of \cite[Lem.~14.8 and Lem.~14.10]{ISem27}. In particular, the argument depends only on $d,\varepsilon_0$ and the constants from \ref{item1: Tb test fct}--\ref{item3: Tb test fct}, and therefore applies verbatim. 
\end{proof}

\begin{proof}[Proof of Corollary~\ref{cor: Kato2}]
    The claim follows directly from Theorem~\ref{thm: Kato property} and Proposition~\ref{prop: homogeneous norms}.
\end{proof}

\section{Applications and consequences of the Kato property}
\label{sec: application}

In this section, we present different applications of the Kato square root property for our fractional operator. These include $\L^2$-boundedness of the Riesz transform associated to $L$ (Section~\ref{sec: Riesz transform}) and topological mapping properties of $L^\frac{1}{2}$ with respect to the coefficients $\mu$ (Section~\ref{sec: holom dependence}) and the fractional parameter $\alpha$ (Section~\ref{sec: alpha dependence}). In the latter, we will also study the behavior of the square root as $\alpha \to 1$ and recover the Kato square root property for elliptic operators $-\div(\mu\nabla\cdot)$ from \cite{AHLMcT}.

\subsection{The Riesz transform associated to \texorpdfstring{$L$}{L}}
\label{sec: Riesz transform}

As an immediate consequence of Theorem~\ref{thm: Kato property}, we have the following properties of the square root.

\begin{lemma}
\label{lem: injectivitiy of sqrt}
    Let $\alpha\in(0,1)$. The square root $L^\frac{1}{2}:\S^{\alpha,2}(\IR^d)\to \L^2(\IR^d)$ is injective but not surjective.
\end{lemma}

\begin{proof}
    For injectivity, we use Theorem~\ref{thm: Kato property}, Proposition~\ref{prop: fourier symbols} and Plancherel's theorem to observe that if $L^\frac{1}{2}u = 0$ for some $u\in \S^{\alpha,2}(\IR^d)$, then
    \begin{align*}
        0=\|L^\frac{1}{2}u\|_{\L^2} \simeq \|\nabla^\alpha u\|_{\L^2} \simeq \| |\xi|^{\alpha}\cF(u)\|_{\L^2}.
    \end{align*}
    As a consequence $\cF (u) = 0$, which implies $u=0$. Next, assume $L^\frac{1}{2}$ is surjective. Then, the open mapping theorem implies that $L^\frac{1}{2}$ is an isomorphism and
    \begin{align*}
        \|u\|_{\S^{\alpha,2}} \simeq \|L^\frac{1}{2}u\|_{\L^2} \simeq \|\nabla^\alpha u\|_{\L^2}
    \end{align*}
    for all $u\in \S^{\alpha,2}(\IR^d)$. However, for $\lambda >0$ and $u\in \C_c^\infty(\IR^d)\setminus \{0\}$, we obtain
    \begin{align*}
        \lambda^\frac{d}{2}\|u\|_{\L^2} = \|u(\cdot/\lambda)\|_{\L^2} \lesssim \|\nabla^\alpha [u(\cdot/\lambda)]\|_{\L^2} = \lambda^{\frac{d}{2}-\alpha}\|\nabla^\alpha u\|_{\L^2}.
    \end{align*}
    Thus, dividing by $\lambda^\frac{d}{2}$ and letting $\lambda \to \infty$, we conclude $\|u\|_{\L^2} =0$, and thus, $u=0$, which is a contradiction.
\end{proof}

\begin{theorem}[Riesz transform on $\L^2$]
\label{thm: riesz transform}
    Let $\alpha_0\in(0,1)$ and $\alpha\in [\alpha_0,1)$. We can extend the fractional Riesz transform
    \begin{align*}
        \nabla^\alpha L^{-\frac{1}{2}}: \cR(L^\frac{1}{2}) \to \L^2(\IR^d;\IC^d)
    \end{align*}
    to a bounded linear operator on $\L^2$ satisfying
    \begin{align}
    \label{eq: two-sided bdd}
        \|\nabla^\alpha L^{-\frac{1}{2}}u\|_{\L^2}\simeq \|u\|_{\L^2} \qquad (u\in\L^2(\IR^d)),
    \end{align}
    where the implicit constants depend only on $d,\alpha_0$ and $\mu_\bullet, \mu^\bullet$.
\end{theorem}

\begin{proof}
By Lemma~\ref{lem: injectivitiy of sqrt}, we notice that $\nabla^\alpha L^{-\frac{1}{2}}: \cR(L^\frac{1}{2}) \to \L^2(\IR^d;\IC^d)$ is well defined. Moreover, by Theorem~\ref{thm: Kato property}, we observe that for every $u= L^\frac{1}{2}v\in \cR(L^\frac{1}{2})$ the two-sided estimate
\begin{align*}
    \|\nabla^\alpha L^{-\frac{1}{2}}u\|_{\L^2} = \|\nabla^\alpha v \|_{\L^2}\simeq \|L^\frac{1}{2}v\|_{\L^2} = \|u\|_{\L^2}
\end{align*}
holds, where the implicit constants depend only on $d,\alpha_0$ and $\mu_\bullet, \mu^\bullet$. Finally, since
\begin{align*}
    \L^2(\IR^d) =\overline{\cR(L)} \subset \overline{\cR(L^\frac{1}{2})} \subset \L^2(\IR^d),
\end{align*}
see \cite[Prop.~2.1.1]{Haase}, we can extend $\nabla^\alpha L^{-\frac{1}{2}}$ to a bounded operator on $\L^2(\IR^d)$ via density, which satisfies \eqref{eq: two-sided bdd} for every $u\in\L^2(\IR^d)$.  
\end{proof}

\subsection{Holomorphic and Lipschitz dependence}
\label{sec: holom dependence}

A particular application that Kato had in mind was the Lipschitz dependence of the square roots $L^{\frac{1}{2}}$ with respect to the coefficients $\mu$ measured in the $\L^{\infty}$-topology, which is proven via holomorphy. To formulate such properties and statements, let us emphasize the coefficients in the notation of $L$ by writing $L_{\mu}$ for the fractional operator with coefficients $\mu$. 

\begin{proposition}
\label{prop: Holomorphic dependence}
Let $\alpha\in(0,1)$. The set
\begin{align*}
    \cO \coloneqq \{ \mu :\mu \text{ satisfies Assumption~\ref{Ass: Coefficients} for some $\mu_{\bullet} , \mu^{\bullet} > 0$} \}
\end{align*}
is open in the $\L^{\infty}$-topology and the map
    \begin{align*}
        \cO \to \mathcal{L}(\S^{\alpha,2}(\IR^d),\L^2(\IR^d)), \quad \mu \mapsto L_{\mu}^{\frac{1}{2}}
    \end{align*}
    is holomorphic, that is, for every $\mu\in\cO$ and $M\in \L^\infty(\IR^d;\IC^{d \times d})$ there exists $r>0$ such that the map
    \begin{align*}
        \{z\in \IC : |z|<r\} \to \cL(\S^{\alpha,2}(\IR^d),\L^2(\IR^d)), \quad z \mapsto L^{\frac{1}{2}}_{\mu+zM}
    \end{align*}
    is holomorphic.
\end{proposition}

This kind of holomorphy of $\mu \mapsto L_{\mu}^{\frac{1}{2}}$ implies the above-mentioned local Lipschitz property for small perturbations, which reads as follows.

\begin{corollary}
\label{cor: Lipschitz estimate}
Let $\alpha\in(0,1)$ and $\mu_1 \in \cO$. Then there exist $\delta>0$ and a constant $C>0$ depending only on $d,\alpha$, $(\mu_1)_\bullet, (\mu_1)^\bullet$ and $\delta$ such that for all $\mu_2\in \cO$ with $\| \mu_2 - \mu_1 \|_{\L^{\infty} (\IR^d ; \IC^{d \times d})}<\delta$ we have
\begin{align*}
 \| L_{\mu_2}^{\frac{1}{2}} - L_{\mu_{1}}^{\frac{1}{2}} \|_{\cL(\S^{\alpha,2} , \L^2)} \leq C \| \mu_2 - \mu_1 \|_{\L^{\infty} (\IR^d ; \IC^{d \times d})}.
\end{align*}
\end{corollary}

We emphasize that this is a nontrivial result even for the class of symmetric coefficients $\mu$. Indeed, already for elliptic operators in divergence form, there is no known proof that establishes the Lipschitz estimate for symmetric coefficients by remaining in the class of operators with symmetric coefficients. All existing proofs use in a crucial way the square root property for operators with nonsymmetric coefficients. 

The proofs of both statements follow verbatim from the corresponding proofs for the generalized Stokes operator (see \cite[Sec.~6]{Haardt_Tolksdorf}), which are based only on the divergence form structure and on the fact that the constants in Theorem~\ref{thm: Kato property} depend on $\mu$ only through $\mu_\bullet$ and $\mu^\bullet$. We therefore omit them and leave the details to the interested reader.

\subsection{Continuous dependence on \texorpdfstring{$\alpha$}{alpha} and the limit \texorpdfstring{$\alpha\to 1$}{alpha to 1}}
\label{sec: alpha dependence}

In this subsection we use the local uniformity in $\alpha$ in Theorem~\ref{thm: Kato property} to recover the Kato square root property for second-order elliptic operators in divergence form in the limit $\alpha\to 1$, and show continuity of the square root with respect to $\alpha$.

Throughout this subsection, $\mu$ is fixed and satisfies Assumption~\ref{Ass: Coefficients}. We write $L_\alpha$ for the operator $L$ from Definition~\ref{def: ses-form and operator} to indicate its dependence on $\alpha$. 

\begin{definition}
Define the sesquilinear form
\begin{align*}
    \a_1 : \W^{1,2}(\IR^d)\times\W^{1,2}(\IR^d)\to\IC,\qquad
    \a_1(u,v)\coloneqq\int\limits_{\IR^d}\mu\nabla u\cdot\overline{\nabla v}\,\d x .
\end{align*}
Additionally, we define the \textit{weak elliptic operator} $\cL_1:\W^{1,2}(\IR^d) \to (\W^{1,2}(\IR^d))'$ by
\begin{align*}
    \langle \cL_1 u, v\rangle \coloneqq \a_1(u,v) \qquad (u,v\in\W^{1,2}(\IR^d))
\end{align*}
and define
\begin{align*}
    \cD(L_1) \coloneqq \{u\in\W^{1,2}(\IR^d): \cL_1u \in \L^2(\IR^d)\}  \quad \text{with} \quad L_1u = \cL_1 u.
\end{align*}
\end{definition}
As in the fractional case, $L_1$ is maximal accretive and sectorial of angle $\omega_0$ from \eqref{def: omega_0} on $\L^2(\IR^d)$. In particular, $L_1$ satisfies the uniform resolvent bounds from Proposition~\ref{prop: L2 resolvent bounds} for $\alpha=1$ and possesses a unique maximal accretive square root $L_1^{1/2}$. We will need the following two convergence results; see also \cite[Prop.~4.4 and Lem.~5.3]{Comi_StefaniI}. We will use the conventions 
\begin{align*}
    \nabla^1 \coloneqq \nabla , \qquad \div^1\coloneqq\div, \qquad \S^{1,2}(\IR^d)\coloneqq \W^{1,2}(\IR^d)
\end{align*}
to unify the notation.

\begin{lemma}
\label{lem: convergence frac ops}
    Let $\alpha\in (0,1]$, $f\in\W^{1,2}(\IR^d)$ and $F\in\C_c^\infty(\IR^d;\IC^d)$. Then, we have
    \begin{align*}
        \nabla^\beta f\to\nabla^\alpha f \quad\text{in }\L^2(\IR^d;\IC^d)
        \qquad\text{and}\qquad
        \div^\beta F\to\div^\alpha F \quad\text{in }\L^2(\IR^d)
    \end{align*}
    as $ \beta\to\alpha$.
\end{lemma}
\begin{proof}
Since $\W^{1,2}(\IR^d)\subset\S^{\alpha,2}(\IR^d)$ for all $\alpha\in(0,1]$ by Lemma~\ref{lem: embedding}, Proposition~\ref{prop: fourier symbols} together with Plancherel's theorem yields
\begin{align*}
    \|\nabla^\beta f-\nabla^\alpha f\|_{\L^2}^2
    =\int\limits_{\IR^d}\big||\xi|^\beta-|\xi|^\alpha\big|^2|\cF(f)(\xi)|^2\,\d\xi .
\end{align*}
The integrand tends to zero as $\beta\to\alpha$ and is bounded by $4(1+|\xi|^2)|\cF(f)(\xi)|^2$, which is integrable, because $f\in\W^{1,2}(\IR^d)$. Hence, the claim follows
from dominated convergence. The same argument applies to the fractional divergence.
\end{proof}

\begin{lemma}
\label{lem: weak-strong resolvent convergence}
    Let $\alpha\in (0,1]$ and $(\alpha_n)_n\subset(0,1]$ with $\alpha_n\to \alpha$, $s>0$ and $f\in\L^2(\IR^d)$. Then
    \begin{align*}
        (1+sL_{\alpha_n})^{-1}f\to (1+sL_\alpha )^{-1}f\quad \text{and} \quad \nabla^{\alpha_n}(1+sL_{\alpha_n})^{-1}f\to\nabla^\alpha(1+sL_\alpha)^{-1}f
    \end{align*}
    in $\L^2$ as $n\to \infty$.
\end{lemma}

\begin{proof}
The proof comes in two steps. The first step shows that the claimed convergence holds weakly in $\L^2$. Subsequently, we upgrade the weak convergence to a strong convergence in $\L^2$ in the second step.\\

\noindent \textbf{Step 1: Weak convergence in $\L^2$.}
    Set $u_n\coloneqq(1+sL_{\alpha_n})^{-1}f$. By Definition~\ref{def: ses-form and operator}, we have $u_n\in\S^{\alpha_n,2}(\IR^d)$ and
    \begin{align}
    \label{eq: resolvent eq n}
        \int\limits_{\IR^d}u_n\overline{v}\,\d x
        +s\int\limits_{\IR^d}\mu\nabla^{\alpha_n}u_n\cdot\overline{\nabla^{\alpha_n}v}\,\d x
        =\int\limits_{\IR^d}f\overline{v}\,\d x
    \end{align}
    for every $v\in\S^{\alpha_n,2}(\IR^d)$.
    Testing with $v=u_n$ and using Assumption~\ref{Ass: Coefficients} yields
    \begin{align*}
        \|u_n\|_{\L^2}^2
        +s\mu_\bullet \|\nabla^{\alpha_n}u_n\|_{\L^2}^2\leq \int\limits_{\IR^d}|u_n|^2\,\d x
        +s\Re\int\limits_{\IR^d}\mu\nabla^{\alpha_n}u_n\cdot\overline{\nabla^{\alpha_n}u_n}\,\d x
        \leq \bigg|\int\limits_{\IR^d}f\overline{u_n}\,\d x\bigg|.
    \end{align*}
    Applying Young's inequality gets us 
    \begin{align*}
        \|u_n\|_{\L^2}^2
        +s\mu_\bullet \|\nabla^{\alpha_n}u_n\|_{\L^2}^2\leq \frac{\|f\|_{\L^2}^2}{2} +\frac{\|u_n\|^2_{\L^2}}{2}
    \end{align*}
    and hence  
    \begin{align}
    \label{eq: uniform bdd u}
        \|u_n\|^2_{\L^2}+ 2s\mu_\bullet \|\nabla^{\alpha_n}u_n\|^2_{\L^2}\leq \|f\|^2_{\L^2}.
    \end{align}
    Thus, after passing to subsequences $(u_{n_k})_k$ and $(\nabla^{\alpha_{n_k}}u_{n_k})_k$, there exist $u\in \L^2(\IR^d)$ and $G\in \L^2(\IR^d;\IC^d)$ such that
    \begin{align}
    \label{eq: weak convergence}
        u_{n_k}\rightharpoonup u \quad \text{and} \quad \nabla^{\alpha_{n_k}}u_{n_k}\rightharpoonup G \qquad\text{weakly in }\L^2.
    \end{align}
    First, we show $u\in\S^{\alpha,2}(\IR^d)$ with $\nabla^\alpha u=G$. For $F\in\C_c^\infty(\IR^d;\IC^d)$, Definition~\ref{def: S space} gives
    \begin{align*}
        \int\limits_{\IR^d}\nabla^{\alpha_{n_k}}u_{n_k}\cdot F\,\d x
        =-\int\limits_{\IR^d}u_{n_k} \, \div^{\alpha_{n_k}}F\,\d x .
    \end{align*}
    Observe that the left-hand side tends to $\int_{\IR^d} G\cdot F\,\d x$ by \eqref{eq: weak convergence}. Moreover, on the right-hand side,
    $u_{n_k}\rightharpoonup u$ weakly in $\L^2$ and $\div^{\alpha_{n_k}} F\to\div^\alpha F$ strongly in $\L^2$ by Lemma~\ref{lem: convergence frac ops}. Thus, it converges to $-\int_{\IR^d} u \, \div^\alpha F\,\d x$. Hence, we obtain
    \begin{align*}
        \int\limits_{\IR^d} G\cdot F\,\d x
        =-\int\limits_{\IR^d}u \, \div^\alpha F\,\d x
    \end{align*}
    for every $F\in\C_c^\infty(\IR^d;\IC^d)$, which shows  $u\in\S^{\alpha ,2}(\IR^d)$ with $\nabla^\alpha u=G$.
    Next, let $v\in\C_c^\infty(\IR^d)$ and split
    \begin{align*}
        \int\limits_{\IR^d}\mu\nabla^{\alpha_{n_k}}u_{n_k}\cdot\overline{\nabla^{\alpha_{n_k}}v}\,\d x
        =\int\limits_{\IR^d}\mu\nabla^{\alpha_{n_k}}u_{n_k}\cdot\overline{\nabla^\alpha  v}\,\d x
        +\int\limits_{\IR^d}\mu\nabla^{\alpha_{n_k}}u_{n_k}\cdot\overline{(\nabla^{\alpha_{n_k}}v-\nabla^\alpha v)}\,\d x .
    \end{align*}
    Notice that the first term on the right-hand side tends to $\int_{\IR^d}\mu\nabla^\alpha  u\cdot\overline{\nabla^\alpha v}\,\d x$ by \eqref{eq: weak convergence}. Furthermore, the second term on the right-hand side is bounded by 
    \begin{align*}
        \bigg|\int\limits_{\IR^d}\mu\nabla^{\alpha_{n_k}}u_{n_k}\cdot\overline{(\nabla^{\alpha_{n_k}}v-\nabla^\alpha v)}\,\d x\bigg|\leq \mu^\bullet \big(\sup\limits_{k\in\IN} \|\nabla^{\alpha_{n_k}}u_{n_k}\|_{\L^2}\big)\|\nabla^{\alpha_{n_k}}v-\nabla^\alpha v\|_{\L^2},
    \end{align*}
    where the supremum is finite, because $(\nabla^{\alpha_{n_k}}u_{n_k})_{k}$ converges weakly in $\L^2$. Hence, it converges to zero as $k\to \infty$ by Lemma~\ref{lem: convergence frac ops}. As a consequence, we obtain
    \begin{align*}
        \int\limits_{\IR^d}\mu\nabla^{\alpha_{n_k}}u_{n_k}\cdot\overline{\nabla^{\alpha_{n_k}}v}\,\d x \to \int\limits_{\IR^d}\mu\nabla^\alpha u\cdot\overline{\nabla^\alpha v}\,\d x
    \end{align*}
    for $k\to \infty$. Finally, passing to the limit in~\eqref{eq: resolvent eq n} yields
    \begin{align*}
        \int\limits_{\IR^d}u\overline{v}\,\d x+s\int\limits_{\IR^d}\mu\nabla^\alpha u\cdot\overline{\nabla^\alpha v}\,\d x=\int\limits_{\IR^d}f\overline{v}\,\d x
    \end{align*}
    for every $v\in\C_c^\infty(\IR^d)$, and by density for all $v\in\S^{\alpha,2}(\IR^d)$. Hence, $u=(1+sL_\alpha)^{-1}f$ and
    \begin{align*}
        (1+sL_{\alpha_{n_k}})^{-1}f \rightharpoonup (1+sL_\alpha)^{-1}f \quad \text{and} \quad \nabla^{\alpha_{n_k}}(1+sL_{\alpha_{n_k}})^{-1}f \rightharpoonup \nabla^\alpha (1+sL_\alpha)^{-1}f 
    \end{align*}
    weakly in $\L^2$. Since every subsequence of $(u_n,\nabla^{\alpha_n}u_n)_n$ satisfies the uniform bound~\eqref{eq: uniform bdd u}, the argument above applies to every subsequence. More precisely, every subsequence of $(u_n,\nabla^{\alpha_n}u_n)_n$ has a further subsequence which converges weakly to $(u,\nabla^\alpha u)$ in $\L^2$, where $u = (1+sL_\alpha)^{-1}f$. As these limits are the same for all subsequences, the whole sequences $(u_n)_n$ and $(\nabla^{\alpha_n}u_n)_n$ has to converge weakly to $u$ and $\nabla^\alpha u$, respectively.\\

    \noindent \textbf{Step 2: Strong convergence in $\L^2$.}
    Testing the resolvent equations for $u_n$ and $u$ with $u_n$ and $u$, respectively, and taking real parts yields
    \begin{align*}
        \|u_n\|_{\L^2}^2+s q(\nabla^{\alpha_n}u_n)=\Re\langle f,u_n\rangle
        \quad\text{and}\quad
        \|u\|_{\L^2}^2+s q(\nabla^{\alpha}u)=\Re\langle f,u\rangle,
    \end{align*}
    where we set
    \begin{align*}
        q(F)\coloneqq\Re\int\limits_{\IR^d}\mu F\cdot\overline{F}\,\d x\qquad (F\in\L^2(\IR^d;\IC^d)).
    \end{align*}
    Subtracting both identities and writing $G_n\coloneqq\nabla^{\alpha_n}u_n$ and $G\coloneqq\nabla^{\alpha}u$, we obtain
    \begin{align*}
        \|u_n\|_{\L^2}^2- \|u\|_{\L^2}^2+ sq(G_n)-sq(G)=\Re\langle f,u_n-u\rangle .
    \end{align*}
    Using the identities
    \begin{align*}
        \|u_n\|_{\L^2}^2- \|u\|_{\L^2}^2 = \|u_n-u\|_{\L^2}^2+2\Re \langle u_n -u, u\rangle
    \end{align*}
    and
    \begin{align*}
        q(G_n)-q(G) = q(G_n-G) + \Re\int\limits_{\IR^d}\mu G\cdot\overline{G_n-G}\,\d x+\Re\int\limits_{\IR^d}\mu (G_n-G)\cdot\overline{G}\,\d x
    \end{align*}
    in the equation above, we arrive after rearranging at
    \begin{align*}
        \|u_n-u\|_{\L^2}^2+ sq(G_n-G)&=\Re\langle f,u_n-u\rangle -2\Re \langle u_n -u, u\rangle\\
        &\quad -s\Re\int\limits_{\IR^d}\mu G\cdot\overline{G_n-G}\,\d x-s\Re\int\limits_{\IR^d}\mu (G_n-G)\cdot\overline{G}\,\d x.
    \end{align*}
    Because $f,u\in\L^2(\IR^d)$ and $G=\nabla^\alpha u \in \L^2(\IR^d;\IC^d)$, the right-hand side of the equality above vanishes for $n\to \infty$ due to the weak convergence from Step~1. Hence, using Assumption~\ref{Ass: Coefficients}, we get
    \begin{align*}
         \|u_n-u\|_{\L^2}^2+s\mu_\bullet\|G_n-G\|_{\L^2}^2\leq \|u_n-u\|_{\L^2}^2+ sq(G_n-G) \to 0
    \end{align*}
    as $n\to \infty$, which proves the claim.
\end{proof}

Combining these two convergence lemmata with the uniform Kato property from Theorem~\ref{thm: Kato property} recovers the classical Kato property for the local elliptic operator $L_1$.

\begin{proof}[Proof of Proposition~\ref{prop: Kato for elliptic}]
As argued in Section~\ref{Sec: Reduction to a square function estimate}, it is well known that the Kato property of $L_1$ follows from the square function estimate
\begin{align}
\label{eq:SFE for elliptic operator}
    \int\limits_0^\infty \|t L_1(1+t^2L_1)^{-1}u\|^2_{\L^2}\,\frac{\d t}{t} &\lesssim \|\nabla u\|_{\L^2}^2,
\end{align}
for all $u\in\W^{1,2}(\IR^d)$, see for instance \cite[Cor.~12.9]{ISem27}. Recall from the proof of Theorem~\ref{thm: Kato property} that for all $\alpha_0\in(0,1)$ and $\alpha\in[\alpha_0,1)$, we proved the square function estimate
\begin{align}
\label{eq: sfe 2}
    \int\limits_0^\infty \|t^\alpha L_\alpha(1+t^{2\alpha}L_\alpha)^{-1}u\|^2_{\L^2}\,\frac{\d t}{t} &\leq C\|\nabla^\alpha u\|_{\L^2}^2 \quad (u\in \S^{\alpha,2}(\IR^d))
\end{align}
for $C>0$ depending only on $d,\alpha_0$ and $\mu_\bullet,\mu^\bullet$.
To this end, fix $\alpha_0\coloneqq\frac{1}{2}$. Let $u\in\W^{1,2}(\IR^d)$ and $(\alpha_n)_n\subset[\frac{1}{2},1)$ with $\alpha_n\to1$. Using $tL_{\alpha_n}(1+t^2L_{\alpha_n})^{-1}=t^{-1}(1-(1+t^2L_{\alpha_n})^{-1})$ and Lemma~\ref{lem: weak-strong resolvent convergence} with $s=t^2$, we have
\begin{align*}
    tL_{\alpha_n}(1+t^2L_{\alpha_n})^{-1}u\to tL_1(1+t^2L_1)^{-1}u
    \qquad \text{in }\L^2(\IR^d).
\end{align*}
Combining this with Fatou's lemma, \eqref{eq: sfe 2} and Lemma~\ref{lem: convergence frac ops}, we obtain
\begin{align*}
    \int\limits_0^\infty\big\|tL_1(1+t^2L_1)^{-1}u\big\|_{\L^2}^2\,\frac{\d t}{t}
    &\leq\liminf_{n\to\infty}\int\limits_0^\infty
    \big\|tL_{\alpha_n}(1+t^2L_{\alpha_n})^{-1}u\big\|_{\L^2}^2\,\frac{\d t}{t}\\
    &=\liminf_{n\to\infty}\alpha_n \int\limits_0^\infty
    \big\|t^{\alpha_n}L_{\alpha_n}(1+t^{2\alpha_n}L_{\alpha_n})^{-1}u\big\|_{\L^2}^2\,\frac{\d t}{t}\\
    &\leq C\liminf_{n\to\infty}\|\nabla^{\alpha_n}u\|_{\L^2}^2\\
    &=C\|\nabla u\|_{\L^2}^2,
\end{align*}
where $C>0$ is uniform with respect to the sequence $(\alpha_n)_{n\in\IN}$. This finishes the proof.
\end{proof}

While the previous statement recovers the square root property at $\alpha=1$ via square function estimates, we can actually prove more. In fact, we show that the square root $L_\alpha^\frac{1}{2}u$ is continuous in $\alpha \in (0,1]$. However, since the domain $\cD(L_\alpha^\frac{1}{2})=\S^{\alpha,2}(\IR^d)$ depends on $\alpha$, we need one fixed common function space. Fortunately, the fractional Sobolev spaces form a continuous embedding chain due to Lemma~\ref{lem: embedding}, so that
\begin{align*}
    \bigcap_{\alpha\in(0,1]}\cD(L_\alpha^\frac{1}{2}) = \W^{1,2}(\IR^d)
\end{align*}
is the natural choice. 

\begin{proof}[Proof of Theorem~\ref{thm: continuity in alpha}]
    Fix $u\in \W^{1,2}(\IR^d)$. Let $\alpha\in(0,1]$ and $(\alpha_n)\subset (0,1]$ with $\alpha_n \to \alpha$. Without loss of generality, we assume $\alpha_n\geq \frac{\alpha}{2}$ for every $n\in\IN$. We divide the proof into two steps.\\
    
    \noindent \textbf{Step 1: Smoothed version.}
    For $\varepsilon>0$, we show the convergence
    \begin{align*}
        L_{\alpha_n}^\frac{1}{2}(1+\varepsilon L_{\alpha_n})^{-1}u \to L_{\alpha}^\frac{1}{2}(1+\varepsilon L_{\alpha})^{-1}u
    \end{align*}
    in $\L^2(\IR^d)$. For this purpose, we use the Balakrishnan representation formula~\cite[Prop.~3.1.12]{Haase}
    \begin{align}
    \label{eq: representation sqrt}
        L_{\alpha_n}^\frac{1}{2}(1+\varepsilon L_{\alpha_n})^{-1}u = \frac{1}{\pi} \int\limits_0^\infty t^{-\frac{1}{2}}(t+L_{\alpha_n})^{-1}L_{\alpha_n}(1+\varepsilon L_{\alpha_n})^{-1}u  \,\d t.
    \end{align}
    Combining the identity $L_{\alpha_n}(s+L_{\alpha_n})^{-1} = 1-s(s+L_{\alpha_n})^{-1}$ for every $s>0$ with the uniform $\L^2$-bounds of Proposition~\ref{prop: L2 resolvent bounds} yields
    \begin{align*}
        \|(t+L_{\alpha_n})^{-1}L_{\alpha_n}(1+\varepsilon L_{\alpha_n})^{-1}u \|_{\L^2} \lesssim \min \Big\{1,\frac{1}{t\varepsilon}\Big\}\|u\|_{\L^2}.
    \end{align*}
    Moreover, by Lemma~\ref{lem: weak-strong resolvent convergence}, the limit
    \begin{align*}
        (t+L_{\alpha_n})^{-1}L_{\alpha_n}(1+\varepsilon L_{\alpha_n})^{-1}u \to  (t+L_{\alpha})^{-1}L_{\alpha}(1+\varepsilon L_{\alpha})^{-1}u
    \end{align*}
    holds in $\L^2(\IR^d)$ as $n\to \infty$. Hence, we can apply dominated convergence to obtain
    \begin{align*}
        L_{\alpha_n}^\frac{1}{2}(1+\varepsilon L_{\alpha_n})^{-1}u &= \frac{1}{\pi} \int\limits_0^\infty t^{-\frac{1}{2}}(t+L_{\alpha_n})^{-1}L_{\alpha_n}(1+\varepsilon L_{\alpha_n})^{-1}u  \,\d t\\
        &\to\frac{1}{\pi} \int\limits_0^\infty t^{-\frac{1}{2}}(t+L_{\alpha})^{-1}L_{\alpha}(1+\varepsilon L_{\alpha})^{-1}u  \,\d t=L_{\alpha}^\frac{1}{2}(1+\varepsilon L_{\alpha})^{-1}u
    \end{align*}
    as $n\to \infty$ for every $\varepsilon>0$.\\

    \noindent \textbf{Step 2: Conclusion.} Because $L_\beta$ is sectorial of angle $\omega_0$ as defined in \eqref{def: omega_0} for all $\beta\in (0,1]$, we have the approximation property
    \begin{align}
    \label{eq: obs 1}
        L_\beta^\frac{1}{2}(1+\varepsilon L_\beta)^{-1}u = (1+\varepsilon L_\beta)^{-1}L_\beta^\frac{1}{2}u  \to L_\beta^\frac{1}{2}u
    \end{align}
    in $\L^2(\IR^d)$ as $\varepsilon\to 0$, see for instance \cite[Prop.~2.1.1]{Haase}. Now, fix some arbitrary $\varepsilon>0$. Using the triangle inequality, we estimate
    \begin{align}
        \label{eq: rom 123}
        \begin{split}
            \|L_\alpha^\frac{1}{2}u - L_{\alpha_n}^\frac{1}{2}u \|_{\L^2} &\leq \|L_\alpha^\frac{1}{2}u - L_{\alpha}^\frac{1}{2}(1+\varepsilon L_\alpha)^{-1}u \|_{\L^2} \\
            &\quad +\|L_{\alpha}^\frac{1}{2}(1+\varepsilon L_\alpha)^{-1}u- L_{\alpha_n}^\frac{1}{2}(1+\varepsilon L_{\alpha_n})^{-1}u \|_{\L^2}\\
            &\quad +\|L_{\alpha_n}^\frac{1}{2}u - L_{\alpha_n}^\frac{1}{2}(1+\varepsilon L_{\alpha_n})^{-1}u \|_{\L^2}\\
            &\eqqcolon \mathrm{I} (\varepsilon)+ \mathrm{II}(n,\varepsilon) + \mathrm{III}(n,\varepsilon).
        \end{split}
    \end{align}
    By Step~1, we know that $\mathrm{II}(n,\varepsilon)\to 0$ as $n\to \infty$.
    Furthermore, by Theorem~\ref{thm: Kato property} with $\alpha_0=\frac{\alpha}{2}$ (or Proposition~\ref{prop: Kato for elliptic} if $\alpha_n=1$ or $\alpha=1$), we have
    \begin{align*}
        \mathrm{III}(n,\varepsilon) &\lesssim \|\nabla^{\alpha_n}(u -(1+\varepsilon L_{\alpha_n})^{-1}u) \|_{\L^2},
    \end{align*}
    where the implicit constant does not depend on $n$ and $\varepsilon$. 
    Thus, applying Lemma~\ref{lem: convergence frac ops} and~\ref{lem: weak-strong resolvent convergence}, we see that 
    \begin{align*}
        \limsup_{n\to\infty}\mathrm{III}(n,\varepsilon) &\lesssim \limsup_{n\to\infty}\|\nabla^{\alpha_n}(u -(1+\varepsilon L_{\alpha_n})^{-1}u) \|_{\L^2}\\
        &= \|\nabla^{\alpha}(u -(1+\varepsilon L_{\alpha})^{-1}u) \|_{\L^2}\\
        &\lesssim \|L_\alpha^\frac{1}{2}u - L_{\alpha}^\frac{1}{2}(1+\varepsilon L_\alpha)^{-1}u \|_{\L^2},
    \end{align*}
    where we used Theorem~\ref{thm: Kato property} (or Proposition~\ref{prop: Kato for elliptic} if $\alpha=1$) in the last step. Using both observations on $\mathrm{II}$ and $\mathrm{III}$ in \eqref{eq: rom 123}, we get
    \begin{align*}
        \limsup_{n\to \infty }\|L_\alpha^\frac{1}{2}u - L_{\alpha_n}^\frac{1}{2}u \|_{\L^2} &\leq \mathrm{I}(\varepsilon)+\limsup_{n\to \infty } \mathrm{III}(n,\varepsilon) \lesssim \mathrm{I}(\varepsilon).
    \end{align*}
    Since this holds for every $\varepsilon>0$, using~\eqref{eq: obs 1} and passing to $\varepsilon\to 0$ yields
    \begin{align*}
        \limsup_{n\to \infty }\|L_\alpha^\frac{1}{2}u - L_{\alpha_n}^\frac{1}{2}u \|_{\L^2} =0,
    \end{align*}
    which proves the claim.
\end{proof}

\medskip

\noindent {\bf AI disclosure statement.}
The author used AI tools for language editing and for checking clarity and internal consistency. The mathematical ideas, proofs and results were developed by the author, who takes full responsibility for the manuscript.

\end{document}